\documentclass[12pt]{article}
\usepackage{subeqn}
\usepackage{graphicx}
\usepackage{ifpdf}
\ifpdf  \fi
\usepackage{pstricks}
\usepackage{amsfonts,amssymb,slashbox}
\usepackage{mathptmx,helvet,courier,makeidx,multicol,footmisc}
\usepackage[numbers]{natbib}
\bibpunct{(}{)}{;}{a}{,}{,}
\usepackage[bookmarksnumbered=true, pdfauthor={Wen-Long Jin}]{hyperref}

\newcommand{\commentout}[1]{}

\newcommand{\ba}{\begin{array}}
        \newcommand{\ea}{\end{array}}
\newcommand{\bc}{\begin{center}}
        \newcommand{\ec}{\end{center}}
\newcommand{\bdm}{\begin{displaymath}}
        \newcommand{\edm}{\end{displaymath}}
\newcommand{\bds} {\begin{description}}
        \newcommand{\eds} {\end{description}}
\newcommand{\ben}{\begin{enumerate}}
        \newcommand{\een}{\end{enumerate}}
\newcommand{\beq}{\begin{equation}}
        \newcommand{\eeq}{\end{equation}}
\newcommand{\bfg} {\begin{figure}[h]}
        \newcommand{\efg} {\end{figure}}
\newcommand{\bi} {\begin {itemize}}
        \newcommand{\ei} {\end {itemize}}
\newcommand{\bqn}{\begin{eqnarray}}
        \newcommand{\eqn}{\end{eqnarray}}
\newcommand{\bqs}{\begin{eqnarray*}}
        \newcommand{\eqs}{\end{eqnarray*}}
\newcommand{\bsl} {\begin{slide}[8.8in,6.7in]}
        \newcommand{\esl} {\end{slide}}
\newcommand{\bss} {\begin{slide*}[9.3in,6.7in]}
        \newcommand{\ess} {\end{slide*}}
\newcommand{\btb} {\begin {table}}
        \newcommand{\etb} {\end {table}}
\newcommand{\bvb}{\begin{verbatim}}
        \newcommand{\evb}{\end{verbatim}}

\newcommand{\m}{\mbox}

\newcommand {\pd}[2] {{\frac {\partial {#1}} {\partial {#2}}}}

\newcommand{\cas}[1]{{{\left \{ \ba #1 \ea \right. }}}

\newcommand{\reff}[1] {{{Figure \ref {#1}}}}
\newcommand{\refe}[1] {{(\ref {#1})}}

\def\a          {{\alpha}}

\def\pmb#1{\setbox0=\hbox{$#1$}%
   \kern-.025em\copy0\kern-\wd0
   \kern.05em\copy0\kern-\wd0
   \kern-.025em\raise.0433em\box0 }

\def\eop{{\hfill $\blacksquare$}}
\def\r{{\rho}}
\newtheorem{theorem}{Theorem}[section]

\newtheorem{corollary}[theorem]{Corollary}

\def\dx     {{\Delta x}}
\def\dt     {{\Delta t}}

\begin{document}
\title{Analysis of kinematic waves arising in diverging traffic flow models}
\author{Wen-Long Jin \footnote{Department of Civil and Environmental Engineering, California Institute for Telecommunications and Information Technology, Institute of Transportation Studies, 4000 Anteater Instruction and Research Bldg, University of California, Irvine, CA 92697-3600. Tel: 949-824-1672. Fax: 949-824-8385. Email: wjin@uci.edu. Corresponding author}}
\maketitle
\begin{abstract}
Diverging junctions are important network bottlenecks, and a better understanding of diverging traffic dynamics has both theoretical and practical implications. In this paper, we first introduce a continuous multi-commodity kinematic wave model of diverging traffic and then present a new framework for constructing kinematic wave solutions to its Riemann problem with jump initial conditions. In supply-demand space, the solutions on a link consist of an interior state and a stationary state, subject to admissible conditions such that there are no positive and negative kinematic waves on the upstream and downstream links respectively. In addition, the solutions have to satisfy entropy conditions consistent with various discrete diverge models. In the proposed analytical framework, kinematic waves on each link can be uniquely determined by the stationary and initial conditions, and we prove that the stationary states and boundary fluxes exist and are unique for the Riemann problem of diverge models when all or partial of vehicles have predefined routes. We show that the two diverge models by Lebacque and Daganzo are asymptotically equivalent. We also prove that the supply-proportional and priority-based diverge models are locally optimal evacuation strategies. With numerical examples, we demonstrate the validity of the analytical solutions of interior states, stationary states, and corresponding kinematic waves. This study presents a unified framework for analyzing traffic dynamics arising in diverging traffic and could be helpful for developing emergency evacuation strategies.
\end{abstract}

\section{{Introduction}}
Essential to effective and efficient transportation control, management, and planning is a better understanding of the evolution of traffic dynamics on a road network, i.e., the formation, propagation, and dissipation of traffic congestion. The seminal work by \citep{lighthill1955lwr,richards1956lwr} (LWR) describes traffic dynamics with kinematic waves, including shock and rarefaction waves, in density ($\r$), speed ($v$), and flux ($q$). Based on a continuous version of traffic conservation, $\pd{\r}{t}+\pd{q}x=0$, and an assumption of a speed-density relationship, $v=V(\r)$, the LWR model can be written as
\bqn
\pd{}{t}\r+\pd {}{x} \r V(\r)=0, \label{lwr}
\eqn
which is for a homogeneous road link with time and location independent traffic characteristics, such as free flow speed, jam density, capacity, and so on. In general, $V(\r)$ is a non-increasing function, and $v_f=V(0)$ is the free flow speed. In addition, $q=Q(\r)\equiv \r V(\r)$ is unimodal with capacity $C=Q(\r_c)$, where $\r_c$ is the critical density. Traffic states with density higher than $\r_c$ are congested or over-critical, and those with lower density are free flowing or under-critical. Here we denote the jam density by $\r_j$, and $\r\in[0,\r_j]$.

In a road network, however, more important and interesting are the formation, propagation, and dissipation of traffic queues caused by network bottlenecks, including merges, diverges, and other network junctions \citep{daganzo1999phase}. But compared with numerous studies on the LWR model and higher-order models of traffic flow on a road link \citep{fhwa2004tft}, studies on traffic dynamics at merging, diverging, and other junctions are scarce. In \citep{fazio1990behavioral}, behavioral models were proposed to capture individual vehicles' diverging maneuvers. In \citep{papageorgiou1990assignment}, diverging flows of vehicles on a path are determined by pre-defined splitting rates. In \citep{daganzo1995ctm}, the First-In-First-Out (FIFO) principle was explicitly introduced so that diverging flows are proportional to turning proportions, which can be time-dependent. But it was noted that the FIFO principle could be violated when one downstream branch is heavily congested. In \citep{liu1996junction,ngoduy2006continuum}, diverging traffic was considered in a so-called friction term of a higher-order model, where diverging flow to an off-ramp is determined by expected diverging flow and the congestion level of the off-ramp. In \citep{munoz2002diverge}, it was shown that First-In-First-Out (FIFO) blockage caused by one congested downstream branch could significantly reduce the discharging flow-rate of the whole diverge, and vehicles may not follow the FIFO principle strictly. Diverging traffic with two or more vehicles have been studied in \citep{daganzo1997special, daganzo1997it, newell1999exit}. In \citep{cassidy2003freeway}, metering strategies were discussed for diverging junctions.
As pointed out in \citep{daganzo1999remarks}, different network bottlenecks can induce different traffic behavior; at diverging junctions, which are different from merging and other junctions, not only capacities of all branches but also the combinations of diverging vehicles on the upstream branch could determine the formation and dissipation of queues. In addition, a better understanding of diverging traffic flow could also lead to more efficient evacuation strategies \citep{Sheffi1982evacuation}.
In this study, we are interested in traffic dynamics arising from diverging junctions for one type of vehicles within the framework of the LWR model.

Considering the analytical power and simplicity of the LWR model, many researchers have attempted to study traffic dynamics arising in general transportation networks in the framework of kinematic wave models. 
In one line, \citet{daganzo1995ctm} and \citet{lebacque1996godunov} extended the Godunov discrete form of the LWR model for computing traffic flows through diverging, diverging, and general junctions. Hereafter we call such models as Cell Transmission Models (CTM). In CTM, so-called traffic demand and supply functions are introduced, and boundary fluxes through various types of junctions can be written as functions of upstream demands and downstream supplies. In CTM, various physically meaningful rules can be used to compute boundary fluxes, such as the First-In-First-Out diverging principle \citep{papageorgiou1990assignment,daganzo1995ctm} and the fair merging principle \citep{jin2003merge}. CTM are discrete in nature and only suitable for numerical simulations. Thus they do not provide any analytical insights on traffic dynamics at a network intersection as the LWR model. 
In another line, \citet{holden1995unidirection} and \citet{coclite2005network} attempted to solve a Riemann problem of an intersection with $m$ upstream links and $n$ downstream links. 
In both of the analytical studies, all links are homogeneous and have the same speed-density relations, and traffic dynamics on each link are described by the LWR model. In \citep{holden1995unidirection}, the Riemann problem with jump initial conditions is solved by introducing an entropy condition that maximizes an objective function of all boundary fluxes. In \citep{coclite2005network}, the Riemann problem is solved to maximize total flux with turning proportions. Both studies were able to describe basic waves arising from a network intersection but also subject to significant shortcomings: (i) All links are assumed to have the same fundamental diagram in both studies; (ii) In \citep{holden1995unidirection}, vehicles can travel to an arbitrary downstream link, and the entropy conditions used are pragmatic and lack of physical interpretations; and (iii) In \citep{coclite2005network}, results are only valid for restricted turning proportions and junctions with no fewer downstream links; i.e., $n\geq m$. In addition, neither of these studies present a unified continuous model of network vehicular traffic.

As in \citep{holden1995unidirection,coclite2005network}, in this study we attempt to analytically obtain kinematic wave solutions of traffic dynamics arising at a diverging junction. However, our study does not bear the same limitations as in these studies: all links can be mainline freeways or off-ramps with the same or different characteristics, and our solutions are physically meaningful and consistent with the discrete supply-demand models of diverging traffic, e.g. those proposed in \citep{daganzo1995ctm,lebacque1996godunov}. 
We first present a continuous kinematic wave model of multi-commodity diverging traffic flow based on the conservation of commodity traffic. 
Following the new framework used to solve Riemann problems for inhomogeneous LWR model at a linear junction \citep{jin2009sd} and for merging traffic flow \citep{jin2010_merge}, we present a new framework for solving the Riemann problem for diverge models. In the Riemann solutions, there can be a stationary state and an interior state for each branch. Here stationary states are the self-similar states at the boundary. That is, in the Riemann solutions, stationary states prevail all links after a long time. In contrast, interior states do not take any space in the continuous solution and only show up in one cell in the numerical solutions as observed in \citep{vanleer1984upwind}. We introduce a so-called supply-demand diagram and discuss the problem in supply-demand space, rather than in $\r-q$ space as in \citep{holden1995unidirection,coclite2005network}. After deriving admissible solutions for upstream and downstream stationary and interior states, we introduce an entropy condition based on various diverge models. We then prove that stationary states and boundary fluxes are unique for given upstream demand and downstream supplies (but interior states may not). 
Then, kinematic waves on a link are determined by the corresponding LWR model with the stationary state and the initial state. 
In a sense, kinematic waves of the Riemann problem can be considered as continuous solutions of the discrete Cell Transmission Model with various diverging rules in \citep{daganzo1995ctm,lebacque1996godunov}.

Different from \citep{holden1995unidirection,coclite2005network}, where the Riemann solutions only comprise of initial and stationary states, here we have additional interior states. Interior states were observed when the inhomogeneous LWR model was used to simulate traffic dynamics on a ring road \citep{jin2003inhlwr,jin2009sd}. Although interior states are not directly related to kinematic waves on all links, they are used in the entropy condition and therefore essential to picking out unique physical solutions. As we can see later, interior states are essential to construct kinematic wave solutions for different diverge models.

The rest of the paper is organized as follows. In Section 2, we introduce a continuous multi-commodity kinematic wave model of diverging traffic. In Section 3, we introduce a new framework for solving the kinematic waves of the Riemann problem with jump initial conditions in supply-demand space. In particular, we derive traffic conservation conditions, admissible conditions of stationary and interior states, and additional entropy conditions based on various discrete diverge models. In Section 4, we solve stationary states and boundary fluxes for diverge models when vehicles have predefined routes. In Section 5, we discuss diverge models in various evacuation strategies. In Section 6, we demonstrate the validity of the proposed analytical framework with numerical examples. In Section 7, we summarize our findings and present some discussions.

\section{{A multi-commodity continuous kinematic wave model of diverging traffic flow and its Riemann problem}}
\bfg\bc
\includegraphics[height=4cm]{diverge20090220figure.22}\caption{An illustration of a diverge network}\label{diverge20090220figure.22}
\ec\efg
We consider a diverge network with $m\geq 2$ downstream links and one upstream link, as shown in \reff{diverge20090220figure.22}. In this network, there are $m+1$ links and $m$ paths. We differentiate all vehicles into $P=m$ commodities according to their paths. We denote the link-commodity incidence variable by $\delta _{p,a} $, which equals 1 if commodity $p$ ($p=1,\cdots,m$) uses link $a$ ($a=1,\cdots,m+1$) and 0 otherwise. Thus $P_a=\sum\nolimits_{p=1}^P {\delta _{p,a} } $ is the number of commodities on link $a$: $P_a=1$ for $a=1,\cdots,m$, and $P_a=m$ for $a=0$.
On a link $a$, the location is denoted by link coordinate $x_a\in[X_a,X_a+L_a]$, where $L_a$ is the length of link $a$, and $x_a=X_a$ and $X_a+L_a$ are the upstream and downstream boundaries respectively.
On the path of a commodity $p$, the location is denoted by commodity coordinate $x_p\in[X_p,X_p+L_p]$, where $L_p=\sum_a \delta_{p,a} L_a$ and we assume that there is no loop on a path. is the length of path $p$, and $x_p=X_p$ and $X_p+L_p$ are at the origin and destination respectively. If $\delta_{p,a}=1$, we denote $L_{p,a}$ as the distance from the origin of path $p$ to the upstream boundary of link $a$, and $x_a$ and $x_p$ follows a one-to-one relation: if $x_p\in[X_p+L_{p,a}, X_p+L_{p,a}+L_a]$, then $x_p$ is on link $a$ and $x_a=x_p-X_p-L_{p,a}+X_a$. That is, $\delta_{p,a}(x_a-X_a-x_p+L_{p,a}+X_p)=0$ for all $a=1,\cdots,m+1$ and $p=1,\cdots,m$

For commodity $p$, we denote density, speed, and flux by $\r_p(x_p, t)$, $v_p(x_p,t)$, and $q_p(x_p,t)=\r_p(x_p, t)v_p(x_p,t)$, respectively. From traffic conservation of commodity $p$, we can have the following continuous conservation equation
\bqn
\pd{\r_p}{t}+\pd{q_p}{x_p}&=&0,\label{mccl}
\eqn
whose derivation is the same as that for single commodity \citep[e.g.][]{haberman1977model,newell1993sim}.
For link $a$, we denote density, speed, and flux by $\r_a(x_a,t)$, $v_a(x_a,t)$, and $q_a(x_a,t)=\r_a(x_a,t) v_a(x_a,t)$,  respectively. Then we have that $\r_a(x_a,t)=\sum_p \delta_{p,a} \r_p(x_a,t)$ and $
q_a(x_a,t)=\sum_p \delta_{p,a} q_p(x_a,t)$. Note that, $\r_p(x_a,t)$ exists only when link $a$ is on path $p$ and  $\r_p(x_a,t)=\r_p(x_p,t)$ with $x_a=x_p-L_{p,a}-X_p+X_a$. It is the same for $v_p(x_a,t)$ and $q_p(x_a,t)$.
We assume that traffic streams of different commodities on link $a$ are homogeneous and share the same speed at the same location and time.  That is, we have the following speed-density relationships \citep{greenshields1935capacity,delcastillo1995fd_empirical}
\bqn v_p(x_a,t)=v_a(x_a,t) =V(x_a,\r_a(x_a,t)).\label{fd}\eqn
Generally, $V_a(x_a, \r_a)$ is non-increasing in $\r_a$, and $Q(x_a,\r_a)\equiv \r_a V(x_a,\r_a)$ is unimodal in $\r_a$ with its maximum as capacity at $x_a$.
We can see that conservation laws of multi-commodity flows in \refe{mccl} lead to the following LWR model
\bqn
 \pd {} {t}\r_a(x_a,t)+\pd{} {x_a}\r_a(x_a,t)V(x_a,\r_a(x_a,t))&=&0,  \label{inhlwr}
 \eqn
which can work for inhomogeneous roads.
Correspondingly, we can have the following traffic conservation equation for commodity $p$  and $x_p\in[X_p+L_{p,a},X_p+L_{p,a}+L_a]$
\bqn
\pd{}{t}\r_p(x_p,t)+\pd{}{x_p}\r_{p}(x_a,t) V(\r_a(x_a,t))&=&0, \quad p=1,\cdots,m \label{mckw}
\eqn
 where $x_a=x_p-X_p-L_{p,a}+X_a$. 
For commodity $p$, the traffic stream evolves on the corresponding path, and we obtain a one-dimensional hyperbolic conservation law. However, all traffic streams interact with each other on the network, and we have a system of network hyperbolic conservation laws. We hereafter call \refe{mckw} as a multi-commodity kinematic wave (MCKW) model of diverging traffic. 

We can see that traffic flow on a road network cannot be modeled by either one-dimensional or two-dimensional conservation laws, since vehicles of different commodities interact with each other on their shared links.
In particular, for a diverge network with $m$ downstream links, traffic streams of $m$ commodities interact with each other on the upstream link. Traffic dynamics inside each link can be studied by the LWR models in \refe{lwr} or \refe{inhlwr}, and the remained task is to study traffic dynamics at the diverging junction. 
Here we consider the Riemann problem for the MCKW model of diverging traffic in \refe{mckw} with jump initial conditions. Without loss of generality, we assume that all links are homogeneous and infinitely long. For link $a=1,\cdots, m+1$, we assume that its flow-density relation is $q_a=Q_a(\r_a)$, critical density $\r_{c,a}$, and its capacity $C_a$.
For the network in \reff{diverge20090220figure.22}, we set $X_p=-\infty$, $X_p+L_p=\infty$, and $x_p=0$ is at the diverging junction for $p=1,\cdots,m$; $X_a=0$ and $X_a+L_a=\infty$ for $a=1,\cdots,m$; and $X_{0}=-\infty$ and $X_{m+1}+L_{m+1}=0$. Therefore, $L_{p,0}=0$, $L_{p,p}=\infty$, and $\delta_{p,a}(x_a-x_p)=0$ for all $a=0$ and $p=1,\cdots,m$. 

For commodity $p=1,\cdots,m$, we have the following jump initial conditions:
\bqn
\r_p(x_p,0)&=&\cas{{ll} \r_{p,L}, & x_p\in(-\infty,0]\\\r_{p,R}, & x_p\in(0,+\infty)}. \label{mckwriemann}
\eqn
Then upstream link $0$ and downstream link $i=1,\cdots,m$  have constant initial conditions:
\bqn
\r_{0}(x_{{0}},0)&=&\r_{{0}}\equiv \sum_p \r_{p,L}, \: x_{{0}}\in(-\infty,0), \label{linkini1}\\
\r_i(x_i,0)&=&\r_{i}\equiv \r_{i,R}, \: x_i\in(0,+\infty), \quad i=1,\cdots,m \label{linkini2} 
\eqn

\section{{An analytical framework}}
For link $a=0,\cdots, m$, we define the following demand and supply functions with all subscript $a$ suppressed \citep{engquist1980calculation,daganzo1995ctm,lebacque1996godunov}
\bqn
D(\r)&=&Q(\min\{\r,\r_c\})=\cas{{ll}Q(\r),&\m{if } \r\leq \r_c\\C,&\m{if }\r\geq \r_c},\nonumber \\&=&\int_0^\r \chi(s) Q'(s) ds=\int_0^\r \max\{Q'(s),0\} ds \label{def:d}\\
S(\r)&=&Q(\max\{\r,\r_c\})=\cas{{ll}Q(\r),&\m{if } \r\geq \r_c\\C,&\m{if }\r\leq \r_c},\nonumber\\&=&C+\int_0^\r (1- \chi(s)) Q'(s) ds=C+\int_0^\r \min\{Q'(s),0\} ds, \label{def:s}
\eqn
where $\chi(\r)$ equals 1 iff $Q'(\r)\geq 0$ and equals 0 otherwise. 

Here we represent a traffic state in supply-demand space as $U=(D,S)$. This is different from many existing studies, in which traffic states are considered in $\r$-$q$ space.
For the demand and supply functions in \refe{def:d} and \refe{def:s}, we can see that $D$ is non-decreasing with $\r$ and $S$ non-increasing.
Thus $D\leq C$, $S\leq C$, $\max \{D,S\}=C$, and flow-rate $q(U)=\min\{D,S\}$. In addition, $D=S=C$ iff traffic is critical; $D<S=C$ iff traffic is strictly under-critical (SUC); $S<D=C$ iff traffic is strictly over-critical (SOC). Therefore, state $U=(D,S)$ is under-critical (UC), iff $S=C$, or equivalently $D\leq S$; State $U=(D,S)$ is over-critical (OC), iff $D=C$, or equivalently $S\leq D$.

In \reff{fig:fd-ds}(b), we draw a supply-demand diagram for the two fundamental diagrams in \reff{fig:fd-ds}(a).
 On the dashed branch of the supply-demand diagram, traffic is UC and $U=(D,C)$ with $D\leq C$; on the solid branch, traffic is OC and $U=(C,S)$ with $S\leq C$. Compared with the fundamental diagram of a road section, the supply-demand diagram only considers its capacity $C$ and criticality, but not other detailed characteristics such as  critical density, jam density, or the shape of the fundamental diagram.
That is, different fundamental diagrams can have the same demand-supply diagram, as long as they have the same capacity and are unimodal, and their critical densities, jam densities, or shapes are not relevant. 
However, given a demand-supply diagram and its corresponding fundamental diagram, the points are one-to-one mapped. 
  
 \bfg\bc $\ba{c@{\hspace{0.3in}}c}
\includegraphics[height=1.8in]{sta20061026figure.2} &
\includegraphics[height=1.8in]{sta20061026figure.1} \\
\multicolumn{1}{c}{\mbox{\bf (a)}} &
    \multicolumn{1}{c}{\mbox{\bf (b)}}
\ea$ \ec \caption{Fundamental diagrams and their corresponding supply-demand diagrams }\label{fig:fd-ds} \efg
 
In supply-demand space, initial conditions in \refe{linkini1} and \refe{linkini2} are equivalent to (Here $i=1,\cdots,m$ if not otherwise mentioned)
\bqn
U_{0}(x_{0},0)&=&(D_{0},S_{0}), \quad x_{0}\in(-\infty, 0),\label{dslinkini1}\\
U_i(x_i,0)&=&(D_i,S_i), \quad x_i\in(0, +\infty). \label{dslinkini2}  \\
\eqn
In the solutions of the Riemann problem for \refe{mckw} with initial conditions (\ref{dslinkini1}-\ref{dslinkini2}), a shock wave or a rarefaction wave could initiate on a link from the diverging junction at $x=0$, and traffic states on all links become stationary after a long time. We hereafter refer to these states as stationary states. At the boundary, there can also exist interior states \citep{vanleer1984upwind,bultelle1998shock}, which take infinitesimal space and only exist in one cell in numerical solutions. We denote the stationary states on upstream link $0$ and downstream link $i$ by $U_0^-$ and $U_{i}^+$, respectively. We denote the interior states on links $0$ and $i$ by $U_0(0^-,t)$ and $U_{i}(0^+,t)$, respectively. The structure of Riemann solutions on upstream and downstream links are shown in \reff{fig:riemannstructure}, where arrows illustrate the directions of possible kinematic waves. 
Then the kinematic wave on upstream link $0$ is the solution of the corresponding LWR model with initial left and right conditions of $U_0$ and $U_0^-$, respectively. Similarly, the kinematic wave on downstream link ${i}$ is the solution of the corresponding LWR model with initial left and right conditions of $U_{i}^+$ and $U_{i}$, respectively. 

Since vehicles' proportions travel forward along vehicles \citep{lebacque1996godunov}, traffic dynamics on the upstream link follow the First-In-First-Out (FIFO) principle \citep{papageorgiou1990assignment}. If the commodity proportions $\xi_i$ are predefined and constant, as a result of the global FIFO principle, in the Riemann solutions we have
\bqn
q_i&=&\xi_i q_0, \label{globalfifo}
\eqn
which serves as the First-In-First-Out principle \citep{papageorgiou1990assignment}.
Also we have that, in the stationary state $U_0^-$, vehicles' proportions are the same as predefined ones. However, we could have different proportions in the interior state $U_0(0^-,t)$ and denote the corresponding proportion of commodity $i$ by $\xi_i(0^-,t)$.

\bfg\bc $\ba{c@{\hspace{0.3in}}c}
\includegraphics[width=2.9in]{diverge20090220figure.4} &
\includegraphics[width=2.9in]{diverge20090220figure.5} \\
\multicolumn{1}{c}{\mbox{\bf (a)}} &
    \multicolumn{1}{c}{\mbox{\bf (b)}}
\ea$ \ec \caption{The structure of Riemann solutions: (a) Upstream link $0$; (b) Downstream link $i$}\label{fig:riemannstructure} \efg

We denote $q_{0\to i}$ as the flux from link $0$ to link $i$ for $t>0$. The fluxes are determined by the stationary states: the out-flux of link $0$ is $q_0=q(U_0^-)$, and the in-flux of link $i$ is $q_{i}=q(U_{i}^+)$. Furthermore, from traffic conservation at a diverging junction, we have at stationary states
\bqn
q_{0\to i}=q_i=q(U_i^+), \quad q_{0}=q(U_{0}^-)=\sum_{i=1}^{m} q(U_i^+). \label{trafficconservation}
\eqn

\subsection{Admissible stationary and interior states}

As observed in \citep{holden1995unidirection,coclite2005network}, the speed of a kinematic wave on an upstream link cannot be positive, and that on a downstream link cannot be negative. We have the following admissible conditions on stationary states.
\begin{theorem} [Admissible stationary states] \label{admissibless} For initial conditions in \refe{dslinkini1} and \refe{dslinkini2},  stationary states are admissible if and only if 
\bqn
U_0^-&=&(D_0,C_0) \m{ or } (C_0,  S_0^-), \label{upstreamss}
\eqn
where $S_0^-<D_0$ , and 
\bqn
U_{i}^+&=&(C_{i},S_{i}) \m{ or } (D_{i}^+, C_{i}), \label{downstreamss}
\eqn
where $D_{i}^+<S_{i}$ .
\end{theorem}
The proof is quite straightforward and omitted here. The regions of admissible upstream stationary states in both supply-demand and fundamental diagrams are shown in \reff{fig:ssupstream}, and the regions of admissible downstream stationary states are shown in \reff{fig:ssdownstream}. From the figures, we can also determine the types and traveling directions of waves with given stationary and initial states on all links. In particular, the types of kinematic waves and the signs of the wave speeds can be determined in the supply-demand diagram, but the absolute values of the wave speeds have to be determined in the fundamental diagram.

\bfg\bc $\ba{c@{\hspace{0.3in}}c}
\includegraphics[height=1.8in]{diverge20090220figure.6} &
\includegraphics[height=1.8in]{diverge20090220figure.8} \\
\multicolumn{1}{c}{\mbox{\bf (a)}} &
    \multicolumn{1}{c}{\mbox{\bf (b)}}	
\ea$ \ec \caption{Admissible stationary states for upstream link $0$: marked by black dots}\label{fig:ssupstream} \efg

\bfg\bc $\ba{c@{\hspace{0.3in}}c}
\includegraphics[height=1.8in]{diverge20090220figure.10} &
\includegraphics[height=1.8in]{diverge20090220figure.12} \\
\multicolumn{1}{c}{\mbox{\bf (a)}} &
    \multicolumn{1}{c}{\mbox{\bf (b)}}
\ea$ \ec \caption{Admissible stationary states for downstream link $i$: marked by black dots}\label{fig:ssdownstream} \efg

{\em Remark 1.} $U_0^-=U_0$ and $U_i^+=U_i$ are always admissible. In this case, the stationary states are the same as the corresponding initial states, and there are no waves.

{\em Remark 2.}
Out-flux $q_0=\min\{D_0^-,S_0^-\}\leq D_0$ and in-flux $q_i=\min\{D_i^+,S_i^+\}\leq S_i$. That is, $D_0$ is the maximum sending flow and $S_i$ is the maximum receiving flow in the sense of \citep{daganzo1994ctm,daganzo1995ctm}.

{\em Remark 3.} In \citep{lebacque2005network}, a so-called ``invariance principle"  is proposed as follows: if $D_0^-=C_0$, then $q(U_0^-)<D_0$; if $S_i^+=C_i$, then $q(U_i^+)<S_0$. We can see that Theorem \ref{admissibless} is consistent with the ``invariance principle".

\begin{corollary} \label{flux2stationary} For the upstream link $0$, $q_0\leq D_0$; $q_0<D_0$ if and only if $U_0^-=(C_0,q_0)$, and $q_0=D_0$ if and only if $U_0^-=(D_0, C_0)$. For the downstream link $i$, $q_i\leq S_i$; $q_i<S_i$ if and only if $U_i^+=(q_i,C_i)$, and $q_i=S_i$ if and only if $U_i^+=(C_i,S_i)$. That is, given out-fluxes and in-fluxes, the stationary states can be uniquely determined.
\end{corollary}

For interior states, the waves of the Riemann problem on link $0$ with left and right initial conditions of $U_0^-$ and  $U_0(0^-,t)$ cannot have negative speeds. Similarly, the waves of the Riemann problem on link $i$ with left and right initial conditions of $U_i(0^+,t)$ and $U_i^+$ cannot have positive speeds. Therefore, interior states $U_0(0^-,t)$ and $U_i(0^+,t)$ should satisfy the following admissible conditions. 
\begin{theorem}[Admissible interior states]\label{intersta}
For asymptotic stationary states $U_0^-$ and $U_i^+$, interior states $U_0(0^-,t)$ and $U_i(0^+,t)$ in \refe{sdentropy} are admissible if and only if  
\bqn
U_0(0^-,t)&=&\cas{{ll}(C_0,S_0^-)=U_0^-, &\m{when }U_0^-\m{ is SOC; i.e., }S_0^-< D_0^-=C_0 \\(D_0(0^-,t),  S_0(0^-,t)), & \m{when }U_0^-\m{ is UC; i.e., }D_0^-\leq S_0^-=C_0} \label{upstreamis}
\eqn
where $S_0(0^-,t) \geq D_0^-$ , and 
\bqn
U_i(0^+,t)&=&\cas{{ll}(D_i^+,C_i)=U_i^+,&\m{when }U_i^+\m{ is SUC; i.e., }D_i^+<S_i^+=C_i \\ (D_i(0^+,t), S_i(0^+,t)), &\m{when }U_i^+\m{ is OC; i.e., }S_i^+\leq D_i^+=C_i } \label{downstreamis}
\eqn
where $D_i(0^+,t)\geq S_i^+$ .
\end{theorem}
The proof is quite straightforward and omitted here. The regions of admissible upstream interior states in both supply-demand and fundamental diagrams are shown in \refe{fig:isupstream}, and the regions of admissible downstream interior states are shown in \refe{fig:isdownstream}. From the figures, we can also determine the types and traveling directions of waves with given stationary and interior states on all links, but these waves are suppressed and cannot be observed, and we are only able to observe possible interior states in numerical solutions.

\bfg\bc $\ba{c@{\hspace{0.3in}}c}
\includegraphics[height=1.8in]{diverge20090220figure.14} &
\includegraphics[height=1.8in]{diverge20090220figure.16} \\
\multicolumn{1}{c}{\mbox{\bf (a)}} &
    \multicolumn{1}{c}{\mbox{\bf (b)}}
\ea$ \ec \caption{Admissible interior states for upstream link $0$: marked by black dots}\label{fig:isupstream} \efg

\bfg\bc $\ba{c@{\hspace{0.3in}}c}
\includegraphics[height=1.8in]{diverge20090220figure.18} &
\includegraphics[height=1.8in]{diverge20090220figure.20} \\
\multicolumn{1}{c}{\mbox{\bf (a)}} &
    \multicolumn{1}{c}{\mbox{\bf (b)}}
\ea$ \ec \caption{Admissible interior states for downstream link $i$: marked by black dots}\label{fig:isdownstream} \efg

{\em Remark 1.} Note that $U_0(0^-,t)=U_0^-$ and $U_{i}(0^+,t)=U_{i}^+$ are always admissible. In this case, the interior states are the same as the stationary states.

\begin{corollary} \label{flux2staint} For upstream link $0$, $q_0\leq D_0$; $q_0<D_0$ if and only if $U_0(0^-,t)=U_0^-=(C_0,q_0)$, and $q_0=D_0$ if and only if $U_0^-=(D_0, C_0)$, and $U_0(0^-,t)=(D_0(0^-,t), S_0(0^-,t))$ with $S_0(0^-,t)\geq D_0$. For the downstream link $i$, $q_i\leq S_i$; $q_i<S_i$ if and only if $U_i(0^+,t)=U_i^+=(q_i,C_i)$, and $q_i=S_i$ if and only if $U_i^+=(C_i,S_i)$, and $U_i(0^+,t)=(D_i(0^+,t),S_i(0^+,t))$ with $D_i(0^+,t) \geq S_i$.
\end{corollary}

\subsection{Entropy conditions consistent with discrete diverge models}
In order to uniquely determine the solutions of stationary states, we introduce a so-called entropy condition in interior states as follows:
\bqn
q_i&=&F_i(U_{0}(0^-,t), U_1(0^+,t),\cdots, U_m(0^+,t),\xi_1(0^-,t), \cdots, \xi_m(0^-,t)).  \label{sdentropy}
\eqn
That is, the entropy condition uses ``local" information in the sense that it determines boundary fluxes from interior states. In the discrete version of \refe{mckw}, the entropy condition is used to determine boundary fluxes from cells contingent to the diverging junction.
Thus, \[F_i(U_{0}(0^-,t), U_1(0^+,t),\cdots, U_m(0^+,t),\xi_1(0^-,t), \cdots, \xi_m(0^-,t))\] in \refe{sdentropy} can be considered as local, discrete flux functions.

In \citep{daganzo1995ctm}, \[F(U_{0}(0^-,t), U_1(0^+,t),\cdots, U_m(0^+,t),\xi_1(0^-,t), \cdots, \xi_m(0^-,t))\] was proposed to solve the following local optimization problem
\bqn
\max_{U_0^-, U_i^+, U_{0}(0^-,t), U_1(0^+,t),\cdots, U_m(0^+,t),\xi_1(0^-,t), \cdots, \xi_m(0^-,t)} \{q_{0} \} \label{optimizationentropy}
\eqn
subject to 
\bqs
&&q_0 \leq  D_0(0^-,t),\\
&&q_{i} \leq  S_{i}(0^+,t),\\
&&\xi_i(0^-,t)=\m{the proportion of vehicles choosing path } i.
\eqs
Thus, we obtain the total flux as
\bqs
F(U_{0}(0^-,t), U_1(0^+,t),\cdots, U_m(0^+,t),\xi_1(0^-,t), \cdots, \xi_m(0^-,t))=\xi_i(0^-,t)\min_{i=1}^m \{ D_0(0^-,t), \frac{S_{i}(0^+,t)}{\xi_i(0^-,t)} \}.
\eqs

In the literature, a number of other diverge models have been proposed. In \citep{lebacque1996godunov}, the upstream demand is split into commodity demands according to predefined turning proportions, and the in-flux of each downstream link is the minimum of its supply and commodity demand. In \citep{jin2003merge}, turning proportions were proposed to be determined by downstream supplies when vehicles have no predefined routes.
In \citep{Sheffi1982evacuation}, turning proportions were proposed to be determined by downstream speeds in a myopic evacuation scheme.
All these local, discrete diverge models can be considered as entropy conditions, so that we have corresponding continuous diverge models  \refe{mckw}.

\subsection{Summary of the solution framework}
To solve the Riemann problem for \refe{mckw} with the initial conditions in \refe{dslinkini1}-\refe{dslinkini2}, we will first find stationary and interior states that satisfy the aforementioned entropy condition, admissible conditions, and traffic conservation equations. Then the kinematic wave on each link will be determined by the Riemann problem of the corresponding LWR model with initial and stationary states as initial conditions. 
Here we will only focus on solving the stationary states on all links, since the kinematic waves of the LWR model have been well studied in the literature.
From all the conditions, we can see that the feasible domains of stationary and interior states are independent of the upstream supply, $S_i$, and the downstream demand, $D_{m+1}$. That is, the same upstream demand and downstream supply will yield the same solutions of stationary and interior states. However, the upstream and downstream wave types and speeds on each can be related to $S_i$ as shown in \reff{fig:ssupstream}(d) and $D_{m+1}$ as shown in \reff{fig:ssdownstream}(d).

\section{Diverge models with predefined turning proportions}
In this paper, we solve the Riemann problem for a diverging junctions with two downstream links; i.e., $m=2$. 
In this section, we consider two entropy conditions, i.e., two diverge model. Here vehicles have predefined routes; i.e., $\xi_i$ are predefined constants, determined by vehicle route choice behaviors.
We attempt to find the relationships between the boundary fluxes and the initial conditions.
\bqn
q_i=\hat F_i(U_0,U_1, U_2).
\eqn
In contrast to local, discrete flux functions $F_i(U_0(0^-,t),U_1(0^-,t), U_2(0^+,t))$, $\hat F_i(U_0,U_1, U_2)$ can be considered as global, continuous.
With the global, continuous fluxes, we can find stationary states from Corollary \ref{flux2stationary}. With the solution framework in the preceding section, we can then find the kinematic waves of the Riemann problem of \refe{mckw} with initial conditions $(U_0, U_1, U_2)$. 

\subsection{Daganzo's diverge model}
In \citep{daganzo1995ctm}, a FIFO diverge model was proposed based on \refe{optimizationentropy}
\bqn
q_0&=&\min\{D_0(0^-,t),  \frac{S_{1}(0^+,t)}{\xi_1(0^-,t)}, \frac{S_{2}(0^+,t)}{\xi_2(0^-,t)} \}, \label{daganzodiverge}
\eqn
and a local FIFO principle
\bqn
q_i&=&\xi_i(0^-,t) q_0.  \label{localfifo}
\eqn
Comparing \refe{localfifo} and \refe{globalfifo}, we obviously have $\xi_i(0^-,t)=\xi_i$. That is, the commodity proportions in the stationary state are the same as predefined.
Thus in Riemann solutions, stationary and interior states have to satisfy \refe{daganzodiverge}, traffic conservation, and the corresponding admissible conditions.

\begin{theorem} \label{thm:daganzodiverge} For the Riemann problem of the MCKW model of merging traffic in \refe{mckw} with initial conditions in \refe{dslinkini1} and \refe{dslinkini2}, boundary fluxes satisfying the entropy condition in \refe{daganzodiverge}, traffic conservation equations, and the corresponding admissible conditions are:
\bqn
q_0=\min\{D_0, \frac{S_1}{\xi_1}, \frac{S_2}{\xi_2}\}, \label{daganzodiverge2}
\eqn
and $q_i=\xi_i q_0$. The corresponding stationary and interior states are in the following:
\ben
\item If $D_0>q_0$, the stationary state of the upstream link is SOC, and $U_0(0^-,t)=U_0^-=(C_0, q_0)$; if $D_0=q_0$, the stationary state of the upstream link is UC, $U_0^-=(D_0,C_0)$, and $U_0(0^-,t)=U_0^-$ or $U_0(0^-,t)=(D_0(0^-,t),S_0(0^-,t))$ with $D_0(0^-,t)> D_0$ and $S_0(0^-,t)\geq D_0$.
\item If $\frac{S_i}{\xi_i}>q_0$, the stationary state of downstream link $i$ ($i=1,2$) is SUC, and $U_i(0^+,t)=U_i^+=(q_i,C_i)$; if $\frac{S_i}{\xi_i}=q_0$, the stationary state of downstream link $i$ is OC, $U_i^+=(C_i,S_i)$, and $U_i(0^+,t)=U_i^+$ or $U_i(0^+,t)=(D_i(0^+,t),S_i(0^+,t))$ with $D_i(0^+,t)\geq S_i$ and $S_i(0^+,t) > S_i$.
\een
\end{theorem}
The proof of the theorem is given in Appendix A. The solutions of fluxes are illustrated in \reff{diverge20090220figure.2}, in which the starting point of an arrow represents the initial condition $(S_1,S_2)$, and the ending point represents the solution $(q_1,q_2)$. That is, in region I, $D_0<\min_i \frac{S_i}{\xi_i}$, and we have $q_i=\xi_i D_0$; in region II, $\frac{C_1}{\xi_1}<\min\{D_0,\frac{C_2}{\xi_2}\}$, and we have $q_1=S_1$ and $q_2=S_1 \frac{\xi_2}{\xi_1}$; in region III, $\frac{C_2}{\xi_2}<\min\{D_0,\frac{C_1}{\xi_1}\}$, and we have $q_2=S_2$ and $q_1=S_2 \frac{\xi_1}{\xi_2}$; on the boundary line between regions I and II, or the boundary line between regions I and III, $q_i=\xi_i D_0$; on the boundary line between regions II and III, $q_i=S_i$. We can see that, in region I, $D_0<S_1+S_2$, and $q_0=\min\{D_0, S_1+S_2\}$. In regions II and III, $q_0<\min\{D_0, S_1+S_2\}$. That is, due to vehicles' route choice behaviors, the capacity of the diverge, $\min\{D_0, S_1+S_2\}$, is generally under-utilized.

\bfg
\bc\includegraphics[height=8cm]{diverge20090220figure.2}\ec
\caption{The solutions of fluxes for a FIFO diverging junction}\label{diverge20090220figure.2}
\efg

Comparing \refe{daganzodiverge} and \refe{daganzodiverge2}, we can see that the global, continuous fluxes have the same functional form as the local, discrete fluxes. In this sense, the FIFO diverge model \refe{daganzodiverge} is ``invariant". Hereafter, we consider a model invariant if and only if the global, continuous fluxes have the same functional form as the local, discrete fluxes.

\subsection{Lebacque's diverge model}\label{section_lebacque}
In \citep{lebacque1996godunov}, the following diverge model was proposed
\bqn
q_i&=&\min\{\xi_i(0^-,t) D_0(0^-,t), S_i(0^+,t) \}, \label{lebacquediverge}
\eqn
and $q_0=q_1+q_2$. 
Compared with Daganzo's model \refe{daganzodiverge}, this model is locally non-FIFO, and its solutions are the following.

\begin{theorem}\label{thm:lebacquediverge} For the Riemann problem of the MCKW model of merging traffic in \refe{mckw} with initial conditions in \refe{dslinkini1} and \refe{dslinkini2}, boundary fluxes satisfying the entropy condition in \refe{lebacquediverge}, traffic conservation equations, and the corresponding admissible conditions are the same as in \refe{daganzodiverge2}; i.e., 
\bqs
q_0=\min\{D_0, \frac{S_1}{\xi_1}, \frac{S_2}{\xi_2}\},
\eqs
and $q_i=\xi_i q_0$. The corresponding stationary and interior states are in the following:
\ben
\item If $D_0>q_0$, the stationary state of the upstream link is SOC, and $U_0(0^-,t)=U_0^-=(C_0, q_0)$; if $D_0=q_0$, the stationary state of the upstream link is UC, $U_0^-=(D_0,C_0)$, and $U_0(0^-,t)=U_0^-$ or $U_0(0^-,t)=(D_0(0^-,t),S_0(0^-,t))$ with $D_0(0^-,t)> D_0$ and $S_0(0^-,t)\geq D_0$.
\item If $\frac{S_i}{\xi_i}>q_0$, the stationary state of downstream link $i$ ($i=1,2$) is SUC, and $U_i(0^+,t)=U_i^+=(q_i,C_i)$; if $\frac{S_i}{\xi_i}=q_0$, the stationary state of downstream link $i$ is OC, $U_i^+=(C_i,S_i)$, and $U_i(0^+,t)=U_i^+$ or $U_i(0^+,t)=(D_i(0^+,t),S_i(0^+,t))$ with $D_i(0^+,t)\geq S_i$ and $S_i(0^+,t) > S_i$.
\item The interior turning proportions $\xi_i(0^-,t)$ can be determined by interior states and stationary states.
\een
\end{theorem}
The proof of the theorem is in Appendix B. The solutions of boundary fluxes can also be illustrated by \reff{diverge20090220figure.2}.

Compared with Daganzo's local FIFO model \refe{daganzodiverge}, Lebacque's diverge model \refe{lebacquediverge} does not satisfy the local optimization condition in \refe{optimizationentropy} or the local FIFO principle \refe{localfifo} and has interior commodity proportions different from the predefined ones. 
Due to the different functional forms of \refe{lebacquediverge} and \refe{daganzodiverge2}, Lebacque's diverge model is not ``invariant". 
However, Lebacque's model has exact the same fluxes, stationary states, and therefore kinematic waves as Daganzo's. That is, \refe{lebacquediverge} converges to \refe{daganzodiverge} asymptotically and continuously and yields globally optimal and FIFO solutions.
In this sense, both models are ``equivalent" globally and continuously.

\section{Diverge models for emergency evacuation}

\subsection{A supply-proportional evacuation strategy}
We consider the a diverging rule proposed in \citep{jin2003merge}, in which
\bqn
q_i&=&\min\{1,\frac{D_0(0^-,t)}{S_1(0^+,t)+S_2(0^+,t)}\} S_i(0^+,t),\quad i=1,2. \label{evacuationdiverge}
\eqn
In this diverging rule, vehicles do not have predefined routes and belong to the same commodity. This diverge model was applied for emergency evacuation situations in a road network \citep{qiu2008see}.
In this model, we have
\bqs
q_0&=&\min\{D_0(0^-,t), S_1(0^+,t)+S_2(0^+,t)\},
\eqs
and the turning proportions are time-dependent
\bqn
\xi_i=\frac{S_i(0^+,t)}{S_1(0^+,t)+S_2(0^+,t)},\quad i=1,2. \label{evacuationpro}
\eqn

Thus in the Riemann solutions, stationary and interior states have to satisfy \refe{evacuationdiverge}, traffic conservation, and the corresponding admissible conditions.

\begin{theorem} \label{thm:evacuationdiverge} For the Riemann problem of the MCKW model of diverging traffic in \refe{mckw} with initial conditions in \refe{dslinkini1} and \refe{dslinkini2}, stationary and interior states satisfying the entropy condition in \refe{evacuationdiverge}, traffic conservation equations, and the corresponding admissible conditions are the following:
\ben
\item When $S_1+S_2 < D_0$, $U_i^+=U_i(0^+,t)=(C_i,S_i)$ ($i=1,2$) and $U_0^-=U_0(0^-,t)=(C_0,S_1+S_2)$;
\item When $S_1+S_2 = D_0$, $U_i^+=U_i(0^+,t)=(C_i,S_i)$ ($i=1,2$), $U_0^-=(D_0, C_0)$, $U_0(0^-,t)=(D_0, C_0)$ or $(D_0(0^-,t), S_0(0^-,t))$ with $D_0(0^-,t)\geq D_0$ and $S_0(0^-,t)>D_0$ when $D_0<C_0$;
\item When $S_i> \frac{C_i}{C_1+C_2}D_0$ ($i=1,2$), $U_0^-=U_0(0^-,t)=(D_0, C_0)$, and $U^-_i=U_i(0^+,t)=(\frac {C_i}{C_1+C_2} D_0, C_i )$.
\item When $S_1+S_2 > D_0$ and $S_i \leq \frac{C_i}{C_1+C_2}D_0$ ($i,j=1$ or 2 and $i\neq j$), $U_0^-=U_0(0^-,t)=(D_0, C_0)$, $U^+_i=(C_i, S_i)$, $U_i(0^+,t)=(C_i,\frac{C_j}{D_0-S_i} S_i)$, and $U^+_j=U_j(0^+,t)=(D_0-S_i,C_j)$.
\een
\end{theorem}

The proof of the theorem is given in Appendix C. 

\begin{corollary} \label{fairflux}
For the Riemann problem of the MCKW model of diverging traffic in \refe{mckw} with initial conditions in \refe{dslinkini1} and \refe{dslinkini2}, boundary fluxes satisfying the entropy condition in \refe{evacuationdiverge}, traffic conservation equations, and the corresponding admissible conditions are the following:
\ben
\item When $S_1+S_2 \leq D_0$, $q_i=S_i$ ($i=1,2$) and $q_0=S_1+S_2$;
\item When $S_i > \frac{C_i}{C_1+C_2}D_0$ ($i=1,2$), $q_i=\frac{C_i}{C_1+C_2}D_0$ and $q_0=D_0$;
\item When $S_1+S_2 > D_0$ and $S_i \leq \frac{C_i}{C_1+C_2}D_0$ ($i,j=1$ or 2 and $i\neq j$),  $q_i=S_i$, $q_j=D_0-S_i$, and $q_0=D_0$.
\een
That is, for $i,j=1$ or 2 and $i\neq j$, $q_0=\min\{D_0, S_1+S_2\}$.
\bqn
q_i&=&\min\{S_i, \max\{ D_0-S_j,\frac{D_0}{C_1+C_2} C_i\}\}. \label{evacuationdiverge2}
\eqn
\end{corollary}
The solutions of fluxes in four different regions are shown in \reff{diverge20090220figure.1}, in which the starting points of arrows represent the initial conditions in $(D_1,D_2)$, and the ending points represent the solutions of fluxes $(q_1,q_2)$. 
In the figure, we can see four regimes: In regime I, both downstream links have OC stationary states; In regime II and IV, one downstream link has SUC and the other OC stationary states; In regime III, both downstream links have SUC stationary states.
Comparing \refe{evacuationdiverge2} and \refe{evacuationdiverge}, we can see that the evacuation diverge model is not ``invariant". Compared with the diverge model in the preceding section, this model is optimal, since $q_0=\min\{D_0,S_1+S_2\}$.

\bfg
\bc\includegraphics[height=8cm]{diverge20090220figure.1}\ec
\caption{Solutions of fluxes for a supply-proportional emergency evacuation diverge model}\label{diverge20090220figure.1}
\efg

\begin{corollary}
If $U_i$ ($i=1,2$) and $U_0$ satisfy
\bqs
\min\{D_i,S_i\}&=&\min\{S_i, \max\{ D_0-S_j,\frac{D_0}{C_1+C_2} C_i\}\},\\
\min\{D_0,S_0\}&=&\min\{S_1+S_2,D_0\},
\eqs
then the unique stationary states are the same as the initial states, and traffic dynamics at the diverging junction are stationary.
\end{corollary}

\subsection{A priority-based evacuation strategy}
Inspired by \refe{evacuationdiverge2}, we propose a priority-based evacuation strategy ($i,j=1,2$ and $i\neq j$)
\bqn
q_i&=&\min\{S_i(0^+,t), \max\{ D_0(0^-,t)-S_j(0^+,t), \alpha_i D_0(0^-,t)\}\},  \label{priorityevacuation}
\eqn
where $\a_i\in[0,1]$ and $\a_1+\a_2=1$. 

\begin{theorem} \label{thm:priorityevacuation}
For the Riemann problem of the MCKW model of diverging traffic in \refe{mckw} with initial conditions in \refe{dslinkini1} and \refe{dslinkini2}, boundary fluxes satisfying the entropy condition in \refe{priorityevacuation}, traffic conservation equations, and the corresponding admissible conditions are given by
\bqn
q_i&=&\min\{S_i, \max\{ D_0-S_j, \alpha_i D_0\}\}, \label{priorityevacuation2}
\eqn
and $q_0=\min\{D_0, S_1+S_2\}$.
\end{theorem}
The proof of the theorem is given in Appendix D. The solutions of fluxes $(q_1,q_2)$ from $(S_1,S_2,D_0)$ are illustrated in \reff{diverge20090220figure.3}. Clearly we can see that fluxes in \refe{evacuationdiverge2} can be considered as a special case when $\xi_i=C_i/(C_1+C_2)$, and the priority-based evacuation diverge model is invariant.
An extreme case is to give one downstream link an absolute priority for evacuation, e.g., $\a_1=1$ and $\a_2=0$. This can happen when link 1 is shorter or less congestion prone. In this case the fluxes in \refe{priorityevacuation2} become
\bqs
q_1&=&\min\{S_1,D_0\},\\
q_2&=&\min\{S_2, \max\{D_0-S_1\}\}.
\eqs

\bfg
\bc\includegraphics[height=8cm]{diverge20090220figure.3}\ec
\caption{Solutions of fluxes for a priority-based evacuation diverge model}\label{diverge20090220figure.3}
\efg

\subsection{A partial evacuation strategy}
By a partial evacuation scenario, we mean that some vehicles have predefined routes and others do not. For example, $\xi_1\in[0,1]$ and $\xi_2\in[0,1]$ are the predefined portions of vehicles choosing link 1 and 2, respectively, but  $\xi_1+\xi_2$ may be smaller than 1. That is, the remaining portion $1-\xi_1-\xi_2$ can take either route. For this scenario, we propose the following evacuation strategy ($i,j=1,2$ and $i\neq j$)
\bqn
q_i&=&\min\{S_i(0^+,t), \frac{1}{\xi_j} S_j(0^+,t)- S_j(0^+,t), \max\{ D_0(0^-,t)-S_j(0^+,t), \alpha_i D_0(0^-,t)\} \}, \label{partialevacuation}
\eqn
where $\a_i\in[\xi_i,1-\xi_j]$ and $\a_1+\a_2=1$.

\begin{theorem} \label{thm:partialevacuation}
For the Riemann problem of the MCKW model of diverging traffic in \refe{mckw} with initial conditions in \refe{dslinkini1} and \refe{dslinkini2},  boundary fluxes satisfying the entropy condition in \refe{partialevacuation}, traffic conservation equations, and the corresponding admissible conditions are given by
\bqn
q_i&=&\min\{S_i, \frac{1}{\xi_j} S_j- S_j, \max\{ D_0-S_j, \alpha_i D_0\} \}, \label{partialevacuation2}
\eqn
and $q_0=q_1+q_2$.
\end{theorem}
The proof of the theorem is omitted here. The solutions of $(q_1,q_2)$ are illustrated in \reff{diverge20090220figure.23}, in which there are six regimes. Furthermore, we can show that (1) $q_i\geq \xi_i q_0$; (2) When $\xi_1+\xi_2=1$; i.e., when all vehicles have predefined routes, \refe{partialevacuation2} is equivalent to \refe{daganzodiverge2}; (3) When $\xi_1=\xi_2=0$; i.e., when all vehicles have no predefined routes, \refe{partialevacuation2} is equivalent to \refe{priorityevacuation2}.
Therefore, \refe{partialevacuation} can be considered as a generalized diverge model, which encapsulates both normal and evacuation diverge model. In addition, the generalized diverge model is invariant.

\bfg
\bc\includegraphics[height=8cm]{diverge20090220figure.23}\ec
\caption{Solutions of fluxes for a generalized diverge model}\label{diverge20090220figure.23}
\efg

\section{Numerical examples}
In this section, we numerically solve various diverge model and demonstrate the validity of our analytical results. 
Here, both links 0 and 1 are two-lane mainline freeways with a corresponding normalized maximum sensitivity fundamental diagram \citep{delcastillo1995fd_empirical} is ($\r\in[0,2]$)
\bqs
Q(\r)&=& \r \left\{ 1-\exp\left[1-\exp\left( \frac14(\frac{2}{\r}-1)\right)\right]\right\}.
\eqs
Link 2 is a one-lane off-ramp with a fundamental diagram as ($\r\in[0,1]$)
\bqs
Q(\r)&=& \frac 12 \r \left\{ 1-\exp\left[1-\exp\left( \frac14(\frac{1}{\r}-1)\right)\right]\right\}.
\eqs
Note that here the free flow speed on the off-ramp is half of that on the mainline freeway, which is 1. Thus we have the capacities $C_0=C_1=4C_2=0.3365$ and the corresponding critical densities $\r_{c0}=\r_{c1}=2\r_{c2}=0.4876$. The length of all three links is the same as $L=10$, and the simulation time duration is $T=360$. Note that here all quantities in this section are normalized and therefore have no units.

In the numerical examples, we discretize each link into $M$ cells and divide the simulation time duration $T$ into $N$ steps. The time step $\dt=T/N$ and the cell size $\dx=L/M$, with $\dt=0.9 \dx$, satisfy the CFL condition \citep{courant1928CFL}
\bqs
v_f \frac{\dt}{\dx}=\frac{\dt}{\dx} =0.9 \leq 1.
\eqs
Then we use the following finite difference equation  for link $i=0,1,2$:
\bqs
\r_{i,m}^{n+1}&=&\r_{i,m}^n+\frac{\dt}{\dx}(q_{i,m-1/2}^n-q_{i,m+1/2}^n),
\eqs 
where $\r_{i,m}^n$ is the average density in cell $m$ of link $i$ at time step $n$, and the boundary fluxes $q_{i,m-1/2}^n$ are determined by supply-demand methods. For example, for downstream links $i=$1 and 2, the out-fluxes are
\bqs
q_{i,m+1/2}^n&=&\min\{D_{i,m}^n,S_{i,m+1}^n\}, \quad m=1,\cdots, M,
\eqs
where $D_{i,m}^n$ is the demand of cell $m$ on link $i$, $S_{i,m+1}^n$ is the supply of cell $m$, and $S_{i,M+1}^n$ is the supply of commodity $i$. 
For link 0, the in-fluxes are
\bqs
q_{0,m-1/2}^n&=&\min\{D_{3,m-1}^n,S_{3,m}^n\}, \quad m=1,\cdots, M,
\eqs
where $D_{0,0}^n$ is the demand at the origin.
Then the in-fluxes of the downstream links and the out-flux of the downstream link are determined by diverge models, which are discrete versions of \refe{sdentropy}:
\bqs
q_{i,1/2}^n&=&F_i(D_{0,M}^n, S_{1,1}^n, S_{2,1}^n),\\
q_{0,M+1/2}^n&=&q_{1,1/2}^n+q_{2,1/2}^n.
\eqs
We also track the commodity proportions in cell $m$ of link 0, $\xi_{i,m}^n$, as follows \citep{jin2004network}
\bqs
\xi_{i,m}^{n+1}&=&\frac{\r_{0,m}^n}{\r_{0,m}^{n+1}} \xi_{i,m}^n +\frac{\dt}{\dx}\frac{q_{0,m-1/2}^n \xi_{i,m-1}^n - q_{0,m+1/2} \xi_{i,m}^n}{\r_{0,m}^{n+1}}.
\eqs
Note that $\xi_i$ is the predefined proportion of commodity $i$.

In our numerical studies, we only consider Lebacque's diverge model \refe{lebacquediverge} and its invariant counterpart \refe{daganzodiverge}. For Lebacque's diverge model, we have
\bqs
q_{i,1/2}^n&=&\min\{D_{0,M}^n \xi_{i,M}^n, S_{i,1}^n\},\\
q_{0,M+1/2}^n&=&q_{1,1/2}^n+q_{2,1/2}^n.
\eqs
In the invariant Daganzo's diverge model, we have ($i,j=1, 2$ and $i\neq j$)
\bqs
q_{0,M+1/2}^n&=&\min\{D_{0,M}^n, \frac{S_{1,1}^n}{\xi_{1,M}^n}, \frac{S_{2,1}^n}{\xi_{2,M}^n}\},\\
q_{i,1/2}^n&=&\xi_{i,M}^n q_{0,M+1/2}^n.
\eqs

\subsection{Kinematic waves, stationary states, and interior states in Lebacque's diverge model}
In this subsection, we study numerical solutions of Lebacque's diverge model in \refe{lebacquediverge}. Initially, links 0 and 1 carry OC flows with $\r_1=\r_3=1$, and 30\% of the vehicles on link 0 diverge to link 2 starting at $t=0$; i.e., $\xi_1=0.7$, and $\xi_2=0.3$. The initial density on link 2 is $\r_2=0.1$. That is, the initial conditions in supply-demand space is 
$U_0=U_1=(0.3365, 0.2473)$ and $U_2=(0.0500, 0.0841)$. 
Here we use the Neumann boundary condition in supply and demand \citep{collela2004_pde}: $D_{0,0}^n=D_{0,1}^n$, $S_{1,M+1}^n=S_{1,M}^n$, and $S_{2,M+1}^n=S_{2,M}^n$. Therefore, we have a Riemann problem here.

In this case, $\frac{S_2}{\xi_2}<D_0<\frac{S_1}{\xi_1}$. Thus according to Theorem \ref{thm:lebacquediverge}, we should have the following stationary and interior states $U_0^-=U_0(0^-,t)=(C_0, \frac{S_2}{\xi_2})$, $U_1^+=U_1(0^+,t)=(\frac{\xi_1}{\xi_2} S_2, C_1)$, and $U_2^+=U_2(0^+,t)=(C_2, S_2)=(C_2,C_2)$. 
 From the LWR model, there should be a back-traveling rarefaction wave on link 0 connecting $U_0$ to $U_0^-$, since $S_0<\frac{S_2}{\xi_2}$; a forward-traveling shock wave on link 1 connecting $U_1^+$ to $U_1$, since $\frac{\xi_1}{\xi_2} S_2< S_1$; and a forward-traveling rarefaction wave on link 2 connecting $U_2^+$ to $U_2$. Furthermore, from \refe{interiorproportion}, we should have that $\xi_1(0^-,t)=0.5833$.

\bfg\bc
\includegraphics[width=4in]{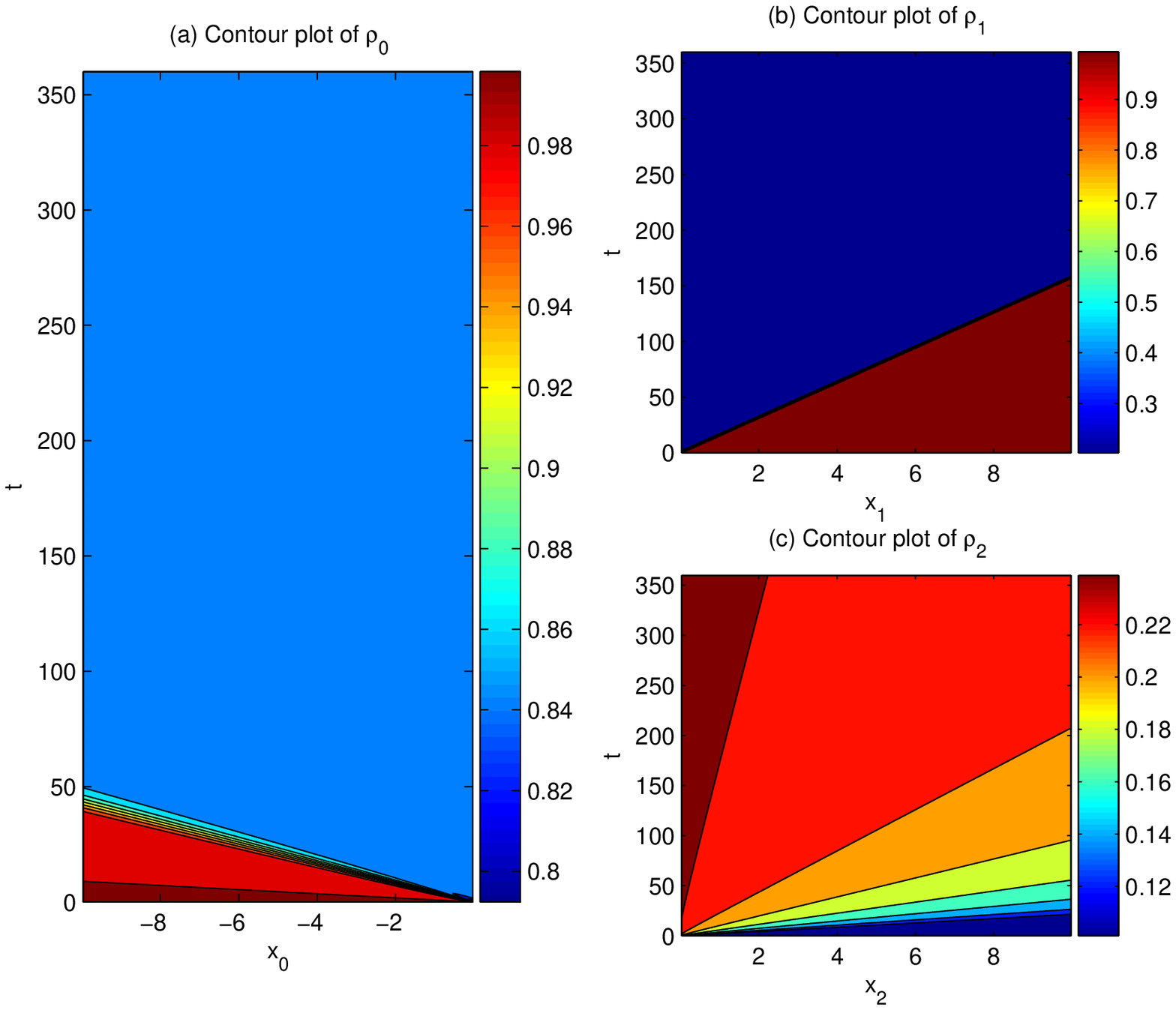} \caption{Solutions of Lebacque's diverge model \refe{lebacquediverge}: $M=160$, $N=6400$.} \label{lebacquediverge20090226contour}
\ec\efg

In \reff{lebacquediverge20090226contour}, the solutions of $\r_0$, $\r_1$, and $\r_2$ are demonstrated with $M=160$ and $N=6400$. From the figures, we can clearly see the predicted kinematic waves. 
In addition, we can observe at $t=T$ the approximate asymptotic values: $U_0^-=U_0(0^-,t)=(0.3365, 0.2804)$, and $\r_0^-=\r_0(0^-,t)=0.8555$; $U_1^+=U_1(0^+,t)=(0.1963, 0.3365)$, and $\r_1^+=\r_1(0^+,t)=0.1963$; and $U_2^+=U_2(0^+,t)=(0.0839,0.0841)$, and $\r_2^+=\r_2(0^+,t)=0.2436\approx \r_{2c}$. These numbers are all very close to the theoretical values and get closer if we reduce $\dx$ or increase $T$. That is, the results are consistent with theoretical results asymptotically.

\bfg\bc
\includegraphics[width=4in]{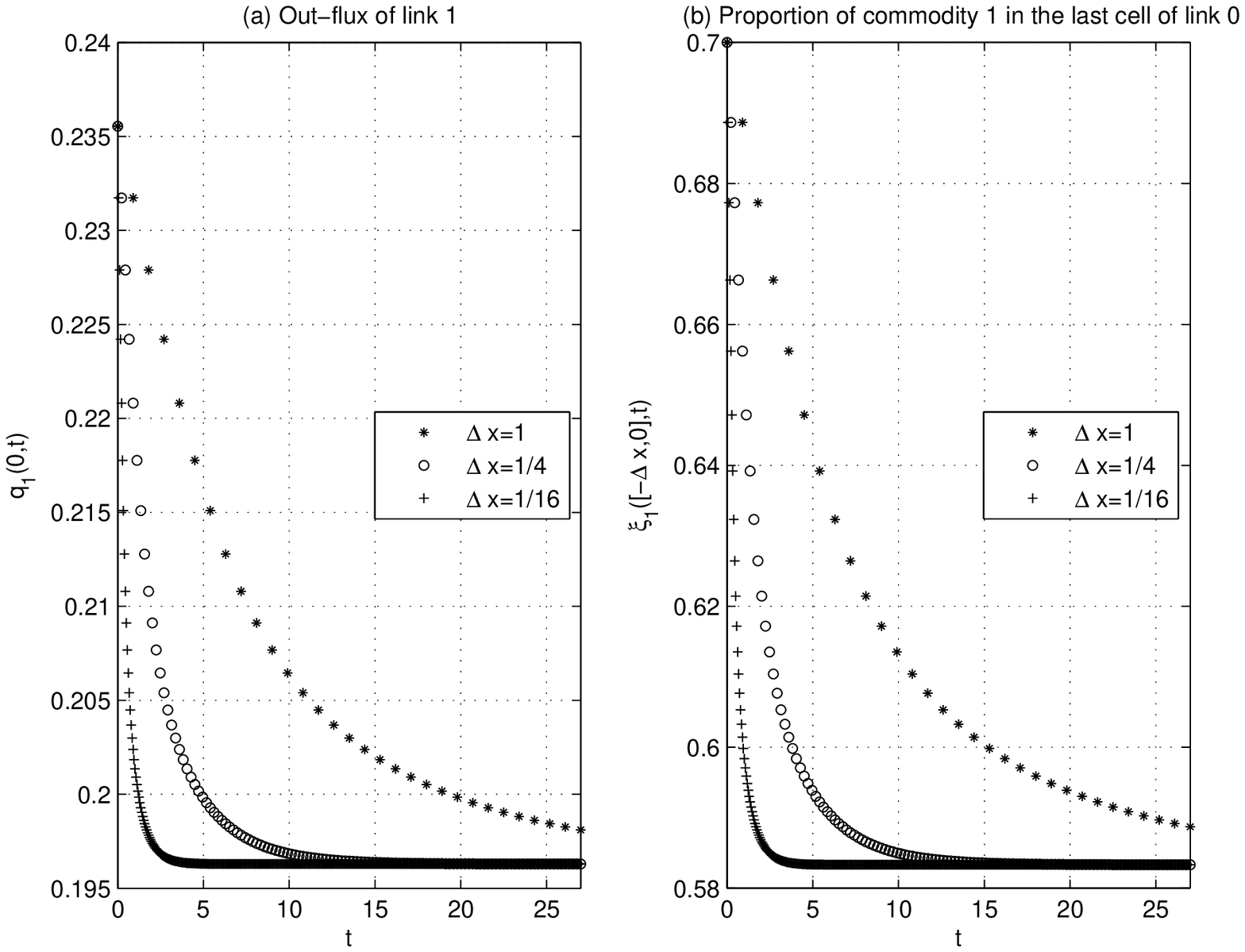} \caption{Evolution of the out-flux and the density in the downstream cell of link 2 for Lebacque's diverge model \refe{lebacquediverge}} \label{lebacquediverge20090226interior}
\ec\efg

In \reff{lebacquediverge20090226interior}, we demonstrate the evolution of the in-flux of link 1 and the proportion of commodity 1 vehicles in the last cell of link 0 for three different cell sizes. From \reff{lebacquediverge20090226interior}(a) we can see that, initially, the out-flux of link 2 is $\min\{\xi_1 D_0, S_1\}=0.2355$, which is not the same but approaches the asymptotic in-flux $\frac{\xi_1}{\xi_2} S_2=0.1963$. Correspondingly the proportion of commodity 1 vehicles in the last cell of link 0 approaches the interior commodity proportion, as shown in \reff{lebacquediverge20090226interior}(b). Similarly, as we decrease the cell size, the numerical results are closer to the theoretical ones at the same time. This figure shows that Lebacque's diverge model is not invariant, but approaches its invariant counterpart asymptotically. Note that the proportion of commodity 1 vehicles in any other cells of link 0 remain constant at 0.7.

\subsection{Comparison of diverge models by Daganzo and Lebacque}\label{compare_daganzo_lebacque}
In this subsection, we compare the numerical solutions of Lebacque's diverge model \refe{lebacquediverge} with its invariant counterpart, Daganzo's diverge model  \refe{daganzodiverge}. Initially, links 0 and 1 carry OC flows with $\r_0=\r_1=1$, and 30\% of the vehicles on link 0 diverge to link 2 starting at $t=0$; i.e., $\xi_1=0.7$, and $\xi_2=0.3$. The initial density on link 2 is $\r_2=0.1$. Different from the example in the preceding subsection, here we use the following boundary conditions: $D_{0,0}^n=D_{0,1}^n$, $S_{1,M+1}^n=S_{1,M}^n$, and $S_{2,M+1}^n=0.05+0.03\sin(n\pi \dt /60)$.  Thus we have a periodic supply on link 2. 

\bfg\bc
\includegraphics[width=4in]{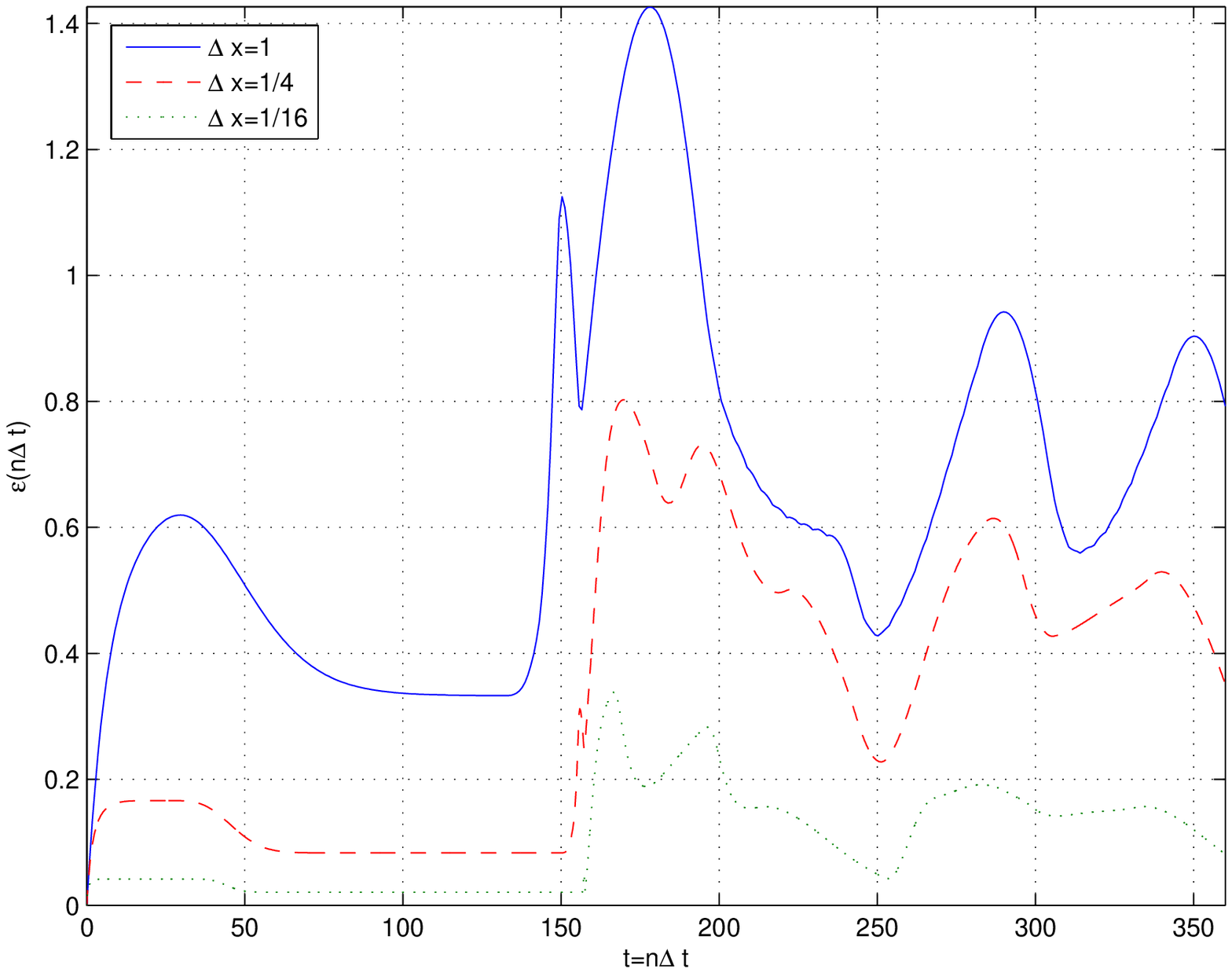}\caption{Difference in the solutions between Lebacque's diverge model \refe{lebacquediverge} and its invariant counterpart  \refe{daganzodiverge}} \label{diverge20090226comparison}
\ec\efg

We use $\r_{i,m}^n$ for the discrete density from Lebacque's diverge model \refe{lebacquediverge} and $\bar \r_{i,m}^n$ from its invariant counterpart  \refe{daganzodiverge}. Then we denote the difference between the two solutions by
\bqn
\epsilon (n \Delta t)&=& \sum_{i=0}^ 2 \sum_{m=1}^M |\r_{i,m}^n- \bar \r_{i,m}^n| \dx.
\eqn
In \reff{diverge20090226comparison}, we can see that the difference decreases if we decreases the cell size. This clearly demonstrates that Lebacque's diverge model \refe{lebacquediverge} converges to its invariant counterpart  \refe{daganzodiverge}.

\section{{Conclusion}}
In this paper, we first introduced a continuous multi-commodity kinematic wave model for a diverge network and defined its Riemann problem. Then, we introduced the supply-demand diagram of traffic flow and proposed a solution framework for the Riemann problem. In the Riemann solutions, each link has two new states: an interior state and a stationary state; and the kinematic waves on a link are determined by the initial state and the stationary state. We then derived admissible conditions for interior and stationary states and introduced entropy conditions consistent with various discrete diverge models. In the analytical framework we proved that the stationary states and boundary fluxes exist and are unique for the Riemann problem for normal diverge, in which vehicles have predefined routes, and evacuation models, in which vehicles may not have predefined routes. With numerical examples, we demonstrated the validity of the solution framework developed here and that Lebacque's diverge model converges to its invariant counterpart, Daganzo's diverge model, when we decrease the cell size.

An important observation is that, for both \refe{lebacquediverge} and \refe{evacuationdiverge}, fluxes computed by discrete supply-demand methods are different from the continuous fluxes. 
For example, the local fluxes from Lebacque's diverge model, \refe{lebacquediverge}, are
\bqs
q_i&=&\min\{\xi_i D_0, S_i\}, \quad i=1,2.
\eqs
When $\xi_i D_0>S_i$; i.e., when the upstream demand is very heavy, we have $q_i=S_i$. In this case, $q_i$ is not proportional to the turning proportion. Thus Lebacque's diverge model violates the FIFO principle. However, from the analysis in Section \ref{section_lebacque} and the numerical example in Section \ref{compare_daganzo_lebacque}, we find that Lebacque's diverge model has the same continuous flux solutions as Daganzo's model, which observes the FIFO principle. Therefore, we conclude that Lebacque's diverge model is not strictly non-FIFO.
As another example, for the supply-proportional evacuation model, at $t=0$ the local fluxes from \refe{evacuationdiverge} are
\bqs
q_i&=&\min\{1,\frac{D_0}{S_1+S_2}\} S_i,\quad i=1,2,
\eqs
which are different from \refe{evacuationdiverge2} when only one downstream is SUC; i.e., when $S_1+S_2>D_0$ and $S_i \leq \frac{C_i}{C_1+C_2}D_0$. However, the analytical results here suggest that the discrete fluxes converge to the continuous ones after a sufficient amount of time or at a given time but with decreasing period of a time interval. 

Comparing kinematic wave solutions of Daganzo's and Lebacque's diverge models, we find that, given the same initial conditions, they have the same stationary states and kinematic wave solutions, but different interior states. In this sense, interior states are essential to distinguish different diverge models. Numerical simulations in Section \ref{compare_daganzo_lebacque} also demonstrate the existence of interior states. Therefore, interior states are essential in understanding diverging traffic flow. This is different from the LWR model for a homogeneous link, in which interior states could exist \citep{jin2003inhlwr,jin2009sd} but are not essential to constructing kinematic wave solutions.

Here we showed that both supply-proportional and priority-based diverge models can be considered locally optimal evacuation strategies. But how to analyze kinematic waves arising in a speed-dependent evacuation model \citep{Sheffi1982evacuation} is subject to further investigations.
In addition to theoretical implications, this study, by improving our understanding of the formation and propagation of traffic congestion caused by diverging bottlenecks, could be helpful for developing, calibrating, and validating diverge models and associated emergency evacuation strategies in the future.	
For example, with different $\a_1$ and $\a_2$, \refe{priorityevacuation} is a priority-based invariant diverge model, which can be used to evacuate vehicles to shorter or less congestion prone links without wasting the capacity of a diverging junction. In the future, we will also be interested in studying kinematic wave solutions of general junctions with multiple upstream and downstream junctions. 

\section*{Acknowledgements}
The author would like to thank two anonymous reviewers for their helpful comments. The views and results contained herein are the author's alone.

\pdfbookmark[1]{Appendix A}{appendixa}
\section*{Appendix A: Proof of Theorem \ref{thm:daganzodiverge}}
{\em Proof}. 
From traffic conservation equations in \refe{trafficconservation}, admissible conditions of stationary states, and the global FIFO principle \refe{globalfifo}, we have $q_i=\xi_i q_0\leq S_i$ and $q_0\leq D_0$. Thus, we have $q_0\leq \min\{D_0, \frac{S_1}{\xi_1}, \frac{S_2}{\xi_2}\}$.

Note that \refe{daganzodiverge} is equivalent to 
\bqs
q_0&=&\min\{D_0(0^-,t),  \frac{S_{1}(0^+,t)}{\xi_1}, \frac{S_{2}(0^+,t)}{\xi_2} \}.
\eqs
We first prove \refe{daganzodiverge2}. 
 Otherwise, $q_0<\min\{D_0, \frac{S_1}{\xi_1}, \frac{S_2}{\xi_2}\}$. (i) Since $q_0<D_0$, from \refe{upstreamss} and \refe{upstreamis}, link $0$ is SOC, and $U_0(0^-,t)=U_0^-=(C_0,q_0)$. (ii) Since $q_i=\xi_i q_0<S_i$, from \refe{downstreamss} and \refe{downstreamis}, link $i$ is SUC, and $U_i(0^+,t)=U_i^+=(\xi_i q_0,C_i)$. Hence from \refe{daganzodiverge} we have $q_0=\min\{C_0,  \frac{C_1}{\xi_1}, \frac{C_2}{\xi_2} \}<D_0\leq C_0$. Thus $q_0=\min\{\frac{C_1}{\xi_1}, \frac{C_2}{\xi_2} \}$. From the FIFO principle \refe{globalfifo} we have
 $q_i=\xi_i q_0=\xi_i(0^-,t) \min\{\frac{C_1}{\xi_1}, \frac{C_2}{\xi_2} \} \geq C_i$, and $q_0=q_i/\xi_i \geq \frac{C_i}{\xi_i}$, which contradicts $q_0\leq \min\{D_0, \frac{S_1}{\xi_1}, \frac{S_2}{\xi_2}\}$.

We consider the following cases.
\bi
\item [(1)] When only one of the three terms on the right hand-side of \refe{daganzodiverge2} equals $q_0$, we have (i) $D_0=q_0<\min\{\frac{S_1}{\xi_1}, \frac{S_2}{\xi_2}\}$, (ii) $\frac{S_1}{\xi_1}=q_0<\min\{D_0, \frac{S_2}{\xi_2}\}$, or (iii) $\frac{S_2}{\xi_2}=q_0<\min\{D_0, \frac{S_1}{\xi_1}\}$. Here we only show the solutions of stationary and interior states for (i), and solutions for (ii) and (iii) can be obtained in a similar fashion. When $\min_i \frac{S_i}{\xi_i}>D_0=q_0$, from \refe{daganzodiverge2} we have $q_0=D_0$, and $q_i=\xi_i q_0<S_i$. From \refe{upstreamss} and \refe{upstreamis}, we have $U_0^-=(D_0,C_0)$, and $U_0(0^-,t)=(D_0(0^-,t),S_0(0^-,t))$ with $S_0(0^-,t) \geq D_0$. From \refe{downstreamss} and \refe{downstreamis}, we have $U_i(0^+,t)=U_i^+=(\xi_i q_0,C_i)$. Then from \refe{daganzodiverge} we have
\bqs
q_0=\min_i\{D_0(0^-,t), \frac{C_i}{\xi_i}\}.
\eqs
Since $\min_i \frac{C_i}{\xi_i}\geq \min_i \frac{S_i}{\xi_i}>D_0=q_0$, we have $D_0(0^-,t)=q_0=D_0$, and $U_0(0^-,t)=U_0^-=(D_0,C_0)$. That is, the upstream and downstream interior states are the same as the corresponding stationary states.
\item [(2)] When two of the three terms on the right hand-side of \refe{daganzodiverge2} equals $q_0$, we have (i) $D_0=\frac{S_1}{\xi_1}=q_0<\frac{S_2}{\xi_2}$, (ii) $D_0=\frac{S_2}{\xi_2}=q_0<\frac{S_1}{\xi_1}$, or (iii) $\frac{S_1}{\xi_1}=\frac{S_2}{\xi_2}=q_0<D_0$. Here we only show the solutions for (i), and solutions for (ii) and (iii) can be obtained in a similar fashion. When $D_0=\frac{S_1}{\xi_1}<\frac{S_2}{\xi_2}$, from Corollary \ref{flux2staint}, we have $U_0^-=(D_0, C_0)$, and $U_0(0^-,t)=(D_0(0^-,t), S_0(0^-,t))$ with $S_0(0^-,t)\geq D_0$; $U_1^+=(C_1,S_1)$, and $U_1(0^+,t)=(D_1(0^+,t),S_1(0^+,t))$ with $D_1(0^+,t) \geq S_1$; and  $U_2(0^+,t)=U_2^+=(q_2,C_2)$. Then from \refe{daganzodiverge} we have
\bqs
q_0=\min\{D_0(0^-,t), \frac{S_1(0^+,t)}{\xi_1},\frac{C_2}{\xi_2}\},
\eqs
which leads to $D_0(0^-,t)=D_0$ or $S_1(0^+,t)=\xi_1 D_0=S_1$. In this case, we can have the following interior states for links 1 and 2: (a) $U_0(0^-,t)=U_0^-$, and $U_1(0^+,t)=(D_1(0^+,t),S_1(0^+,t))$ with $D_1(0^+,t)\geq S_1$ and $S_1(0^+,t) > S_1$; (b) $U_1(0^+,t)=U_1^+$, and $U_0(0^-,t)=(D_0(0^-,t),S_0(0^-,t))$ with $D_0(0^-,t)> D_0$ and $S_0(0^-,t)\geq D_0$.
\item [(3)] When all the three terms on the right hand-side of \refe{daganzodiverge2} equals $q_0$, we have $D_0=\frac{S_1}{\xi_1}=\frac{S_2}{\xi_2}=q_0$. From Corollary \ref{flux2staint}, we have $U_0^-=(D_0, C_0)$, and $U_0(0^-,t)=(D_0(0^-,t), S_0(0^-,t))$ with $S_0(0^-,t)\geq D_0$; $U_1^+=(C_1,S_1)$, and $U_1(0^+,t)=(D_1(0^+,t),S_1(0^+,t))$ with $D_1(0^+,t) \geq S_1$; and $U_2^+=(C_2,S_2)$, and $U_2(0^+,t)=(D_2(0^+,t),S_2(0^+,t))$ with $D_2(0^+,t) \geq S_2$. In this case, at least one of $U_0(0^-,t)=U_0^-$, $U_1(0^+,t)=U_1^+$, and $U_2(0^+,t)=U_2^+$ should be satisfied, and the other can be $U_0(0^-,t)=(D_0(0^-,t),S_0(0^-,t))$ with $D_0(0^-,t)> D_0$ and $S_0(0^-,t)\geq D_0$, $U_1(0^+,t)=(D_1(0^+,t),S_1(0^+,t))$ with $D_1(0^+,t)\geq S_1$ and $S_1(0^+,t) > S_1$, or $U_2(0^+,t)=(D_2(0^+,t),S_2(0^+,t))$ with $D_2(0^+,t)\geq S_2$ and $S_2(0^+,t) > S_2$.
\ei
\eop

\pdfbookmark[1]{Appendix B}{appendixb}
\section*{Appendix B: Proof of Theorem \ref{thm:lebacquediverge}}
{\em Proof}. From traffic conservation equations in \refe{trafficconservation}, admissible conditions of stationary states, and the global FIFO principle \refe{globalfifo}, we have $q_i=\xi_i q_0\leq S_i$ and $q_0\leq D_0$. Thus, we have $q_0\leq \min\{D_0, \frac{S_1}{\xi_1}, \frac{S_2}{\xi_2}\}$. 

We first prove that $q_0$ is given by \refe{daganzodiverge2}. 
 Otherwise, $q_0<\min\{D_0, \frac{S_1}{\xi_1}, \frac{S_2}{\xi_2}\}$, which leads to $q_0<D_0$ and $q_i=\xi_i q_0<S_i$. From Corollary \ref{flux2staint}, we have $U_0(0^-,t)=U_0^-=(C_0,q_0)$ and $U_i(0^+,t)=U_i^+=(q_i,C_i)$.
 Then from \refe{lebacquediverge} we have $q_i=\min\{\xi_i(0^-,t) C_0, C_i \}<S_i\leq S_i$. Thus $q_i=\xi_i(0^-,t) C_0$, and $q_0=q_1+q_2=C_0$, which contradicts $q_0<D_0\leq C_0$.

We consider the following cases.
\bi
\item [(1)] When only one of the three terms on the right hand-side of \refe{daganzodiverge2} equals $q_0$, we have (i) $D_0=q_0<\min\{\frac{S_1}{\xi_1}, \frac{S_2}{\xi_2}\}$, (ii) $\frac{S_2}{\xi_2}=q_0<\min\{D_0, \frac{S_1}{\xi_1}\}$, or (iii) $\frac{S_1}{\xi_1}=q_0<\min\{D_0, \frac{S_2}{\xi_2}\}$. Here we only show the solutions of stationary and interior states for (i) and (ii), and solutions for (iii) can be obtained in a similar fashion. 
\bi
\item [(i)] When $\min_i \frac{S_i}{\xi_i}>D_0=q_0$, from \refe{daganzodiverge2} we have $q_0=D_0$, and $q_i=\xi_i q_0<S_i$. From Corollary \ref{flux2staint}, we have $U_0^-=(D_0,C_0)$, $U_0(0^-,t)=(D_0(0^-,t),S_0(0^-,t))$ with $S_0(0^-,t) \geq D_0$, and $U_i(0^+,t)=U_i^+=(\xi_i q_0,C_i)$. Then from \refe{lebacquediverge} we have
\bqs
q_i=\min\{\xi_i(0^-,t) D_0(0^-,t), C_i\}=\xi_i(0^-,t) D_0(0^-,t) <C_i.
\eqs
Then we have $q_0=D_0(0^-,t)=D_0$, and $U_0(0^-,t)=U_0^-=(D_0,C_0)$. Further we have $\xi_i(0^-,t)=\xi_i$. In this case, the upstream and downstream interior states are the same as the corresponding stationary states.

\item [(ii)] When $\frac{S_2}{\xi_2}<\min\{D_0, \frac{S_1}{\xi_1}\}$, from \refe{daganzodiverge2} we have $q_0=\frac{S_2}{\xi_2}<D_0$, $q_2=S_2$, $q_1=\frac{\xi_1}{\xi_2}S_2<S_1$. From Corollary \ref{flux2staint}, we have $U_0(0^-,t)=U_0^-=(C_0,\frac{S_2}{\xi_2})$, $U_2^+=(C_2,S_2)$, $U_2(0^+,t)=(D_2(0^+,t),S_2(0^+,t))$ with $D_2(0^+,t)\geq S_2$, and $U_1(0^+,t)=U_1^+=(\frac{\xi_1}{\xi_2}S_2, C_1)$. Then from \refe{lebacquediverge} we have
\bqs
q_2&=&\min\{\xi_2(0^-,t) C_0, S_2(0^+,t)\},\\
q_1&=&\min\{\xi_1(0^-,t) C_0, C_1\}=\xi_1(0^-,t) C_0 <C_1.
\eqs
Thus we have
\bqn
\xi_1(0^-,t)&=& \frac{\xi_1 S_2}{\xi_2 C_0}, \label{interiorproportion}
\eqn
and
\bqs
q_0&=&\min\{C_0, S_2(0^+,t)+\xi_1(0^-,t) C_0\}=S_2(0^+,t)+\xi_1(0^-,t) C_0 <C_0.
\eqs
Thus, $q_2=S_2(0^+,t)=S_2$, and $U_2(0^+,t)=U_2^+=(C_2,S_2)$.

\ei
\item [(2)] When two of the three terms on the right hand-side of \refe{daganzodiverge2} equals $q_0$, we have (i) $D_0=\frac{S_1}{\xi_1}=q_0<\frac{S_2}{\xi_2}$, (ii) $D_0=\frac{S_2}{\xi_2}=q_0<\frac{S_1}{\xi_1}$, or (iii) $\frac{S_1}{\xi_1}=\frac{S_2}{\xi_2}=q_0<D_0$. Here we only show the solutions for (i), and solutions for (ii) and (iii) can be obtained in a similar fashion. When $D_0=\frac{S_1}{\xi_1}<\frac{S_2}{\xi_2}$, from Corollary \ref{flux2staint}, we have $U_0^-=(D_0, C_0)$, and $U_0(0^-,t)=(D_0(0^-,t), S_0(0^-,t))$ with $S_0(0^-,t)\geq D_0$; $U_1^+=(C_1,S_1)$, and $U_1(0^+,t)=(D_1(0^+,t),S_1(0^+,t))$ with $D_1(0^+,t) \geq S_1$; and  $U_2(0^+,t)=U_2^+=(q_2,C_2)$. Then from \refe{daganzodiverge2} we have
\bqs
q_1&=&\min\{\xi_1(0^-,t) D_0(0^-,t), S_1(0^+,t)\},\\
q_2&=&\min\{\xi_2(0^-,t) D_0(0^-,t), C_2\}=\xi_2(0^-,t) D_0(0^-,t) <C_2.
\eqs
Thus $D_0(0^-,t)=D_0$ or $S_1(0^+,t)=\xi_1 D_0=S_1$. In this case, we can have the following interior states for links 1 and 2: (a) $U_0(0^-,t)=U_0^-$, $\xi_i(0^-,t)=\xi_i$, and $U_1(0^+,t)=(D_1(0^+,t),S_1(0^+,t))$ with $D_1(0^+,t)\geq S_1$ and $S_1(0^+,t) > S_1$; (b) $U_1(0^+,t)=U_1^+$, $\xi_2(0^-,t)=\xi_2 D_0/D_0(0^-,t)$, and $U_0(0^-,t)=(D_0(0^-,t),S_0(0^-,t))$ with $D_0(0^-,t)> D_0$ and $S_0(0^-,t)\geq D_0$.
\item [(3)] When all the three terms on the right hand-side of \refe{daganzodiverge2} equals $q_0$, we have $D_0=\frac{S_1}{\xi_1}=\frac{S_2}{\xi_2}=q_0$. From Corollary \ref{flux2staint}, we have $U_0^-=(D_0, C_0)$, and $U_0(0^-,t)=(D_0(0^-,t), S_0(0^-,t))$ with $S_0(0^-,t)\geq D_0$; $U_1^+=(C_1,S_1)$, and $U_1(0^+,t)=(D_1(0^+,t),S_1(0^+,t))$ with $D_1(0^+,t) \geq S_1$; and $U_2^+=(C_2,S_2)$, and $U_2(0^+,t)=(D_2(0^+,t),S_2(0^+,t))$ with $D_2(0^+,t) \geq S_2$. In this case, at least one of $U_0(0^-,t)=U_0^-$, $U_1(0^+,t)=U_1^+$, and $U_2(0^+,t)=U_2^+$ should be satisfied, and the other can be $U_0(0^-,t)=(D_0(0^-,t),S_0(0^-,t))$ with $D_0(0^-,t)> D_0$ and $S_0(0^-,t)\geq D_0$, $U_1(0^+,t)=(D_1(0^+,t),S_1(0^+,t))$ with $D_1(0^+,t)\geq S_1$ and $S_1(0^+,t) > S_1$, or $U_2(0^+,t)=(D_2(0^+,t),S_2(0^+,t))$ with $D_2(0^+,t)\geq S_2$ and $S_2(0^+,t) > S_2$. Here $\xi_i(0^-,t)$ can be determined once the interior states are determined.
\ei
\eop

\pdfbookmark[1]{Appendix C}{appendixc}
\section*{Appendix C: Proof of Theorem \ref{thm:evacuationdiverge}}
{\em Proof}. From traffic conservation equations in \refe{trafficconservation} and admissible conditions of stationary states, we can see that 
\bqs
q_0\leq \min\{S_1+S_2, D_0 \}.
\eqs
We first demonstrate that it is not possible that $q_0 < \min\{S_1+S_2, D_0 \} \leq \min\{C_1+C_2, C_0\}$. Otherwise, from \refe{downstreamss} and \refe{downstreamis} we have $U_0(0^-,t)=U_0^-=(C_0,q_0)$ with $q_0<D_0$; Since $q(U_1^+)+q(U_2^+)=q_0<S_1+S_2$, then we have $q(U_i^+)<S_i$ for at least one downstream link, e.g., $q_1<S_1$. From \refe{upstreamss} and \refe{upstreamis} we have $U_1(0^+,t)=U_1^+=(q_1,C_1)$.
Then from the entropy condition in \refe{evacuationdiverge} we have
\bqs
q_0&=&\min\{C_1+S_2(0^+,t), C_0\},\\
q_1&=&\min\{1,\frac{C_0}{C_1+S_2(0^+,t)}\} C_1.
\eqs
Since $q_0<C_0$, from the first equation we have $q_0=C_1+S_2(0^+,t)<C_0$, and from the second equation we have $q_1=C_1$, which contradicts $q_1<S_1$. Therefore, 
\bqs
q_1+q_2=q_0 = \min\{S_1+S_2, D_0 \}.
\eqs
That is, the diverge model \refe{evacuationdiverge} yields the optimal fluxes for any initial conditions.

\bi
\item [(1)] When $S_1+S_2<D_0$, we have $q_0=S_1+S_2<D_0$. We have $U_0(0^-,t)=U_0^-=(S_1+S_2,C_0)$. 
Since $q_1+q_2=S_1+S_2$ and $q_i\leq S_i$, we have $q_i= S_i$, and $U_i^+=(C_i,S_i)$. From \refe{downstreamis} we have $U_i(0^+,t)=(D_i(0^+,t),S_i(0^+,t))$ with $D_i(0^+,t)\geq S_i^+=S_i$. From \refe{evacuationdiverge} we have
\bqs
q_0&=&\min\{S_1(0^+,t)+S_2(0^+,t), C_0\}=S_1+S_2<D_0\leq C_0,\\
q_i&=&\min\{1,\frac{C_0}{S_1(0^+,t)+S_2(0^+,t)}\} S_i(0^+,t) =S_i.
\eqs
Thus, $S_i(0^+,t)=S_i\leq S_i(0^+,t)$. Then $U_i(0^+,t)=U_i^+=(C_i,S_i)$. In this case, there are no interior states on all links.

\item [(2)] When $S_1+S_2=D_0$, we have $q_0=D_0$, and $q_i=S_i$. We have $U_0^-=(D_0,C_0)$ and $U_0(0^-,t)=(D_0(0^-,t),D_0(0^-,t))$ with $S_0(0^-,t)\geq D_0^-=D_0$, and $U_i^+=(C_i,S_i)$ and $U_i(0^+,t)=(S_i(0^+,t),S_i(0^+,t))$ with $D_i(0^+,t)\geq S_i^+=S_i$. 
From \refe{evacuationdiverge} we have
\bqs
q_0&=&\min\{S_1(0^+,t)+S_2(0^+,t), D_0(0^-,t)\}=S_1+S_2=D_0,\\
q_i&=&\min\{1,\frac{D_0(0^-,t)}{S_1(0^+,t)+S_2(0^+,t)}\} S_i(0^+,t) =S_i.
\eqs
We can have the following two scenarios. 
\bi
\item [(2-i)] If $S_1(0^+,t)+S_2(0^+,t)\geq D_0(0^-,t)=S_1+S_2=D_0 \leq S_0(0^-,t)$, then $U_0(0^-,t)=U_0^-=(D_0,C_0)$ and there is no interior state on link 0. Moreover, we have 
\bqs \frac{S_1+S_2}{S_1(0^+,t)+S_2(0^+,t)} S_i(0^+,t) =S_i,
\eqs
 which leads to $S_i(0^+,t)\leq S_i$. From the assumption that $S_1(0^+,t)+S_2(0^+,t)\geq S_1+S_2$, we have $S_i(0^+,t) = S_i$. Further we have $U_i(0^+,t)=U_i^+=(C_i,S_i)$, and there are no interior states on links 1 or 2.
\item [(2-ii)] If $D_0(0^-,t)>S_1(0^+,t)+S_2(0^+,t)=S_1+S_2=D_0$, $S_i(0^+,t) =S_i$. Thus $U_i(0^+,t)=U_i^+=(C_i,S_i)$, and there are no interior states on links 1 or 2. Moreover, $U_0(0^-,t)$ satisfies $S_0(0^-,t)>D_0$ and $D_0(0^-,t)\geq D_0$. Thus there can be multiple interior states on link 0 when $D_0<C_0$.
\ei

\item [(3,4)]
When $S_1+S_2>D_0$, then $q_0=q_1+q_2=D_0$. We have $U_0^-=(D_0,C_0)$, and  $U_0(0^-,t)=(D_0(0^-,t),D_0(0^-,t))$ with $S_0(0^-,t)\geq D_0^-=D_0$. For downstream links, at least one of the stationary states is SUC. Otherwise, from \refe{downstreamss} we have $U_i^+=(C_i,S_i)$, and $q_1+q_2=S_1+S_2>D_0$, which is impossible.  	 
From \refe{evacuationdiverge} we have
\bqs
q_0&=&\min\{S_1(0^+,t)+S_2(0^+,t), D_0(0^-,t)\}=D_0<S_1+S_2,\\
q_i&=&\min\{1,\frac{D_0(0^-,t)}{S_1(0^+,t)+S_2(0^+,t)}\} S_i(0^+,t).
\eqs
If $S_1(0^+,t)+S_2(0^+,t)\leq D_0(0^-,t)$, then $S_1(0^+,t)+S_2(0^+,t)=D_0<S_1+S_2$ and $q_i=S_i(0^+,t)$. This is not possible for the SUC stationary state $U_i^+=U_i(0^+,t)=(q_i,C_i)$ with $q_i<S_i\leq C_i$. Thus $D_0(0^-,t)<S_1(0^+,t)+S_2(0^+,t)$, $D_0(0^-,t)=D_0<S_1+S_2$, and $U_0(0^-,t)=U_0^-=(D_0,C_0)$. Hence for both downstream links
\bqs
q_i&=&\frac{D_0}{S_1(0^+,t)+S_2(0^+,t)}S_i(0^+,t).
\eqs

\bi
\item [(3)] When $S_i> \frac{C_i}{C_1+C_2}D_0$ ($i=1,2$), stationary states on both links 1 and 2 are SUC with $U_i^+=U_i(0^+,t)=(q_i,C_i)$ with $q_i<S_i$. Otherwise, we assume that link 1 is SUC with $U_1(0^+,t)=U_1^+=(q_1,C_1)$ and link 2 is OC with $U_2^+=(S_2,C_2)$. Then 
\bqs
S_2&=&\frac{D_0}{C_1+S_2(0^+,t)}S_2(0^+,t) \leq \frac{D_0}{C_1+C_2}C_2<S_2,
\eqs
which is impossible.
From \refe{evacuationdiverge}, we have 
\bqs
q_i=\frac{D_0}{C_1+C_2} C_i,
\eqs
and $U_i(0^+,t)=U_i^+=(q_i,C_i)$.

\item[(4)]
When $S_1+S_2 > S_{3}$ and $S_i \leq \frac{C_i}{C_1+C_2}D_0$ ($i,j=1$ or 2 and $i\neq j$), we can show that stationary states on links $j$ and $i$ are SUC and OC respectively with $U_j^+=U_j(0^+,t)=(q_j,C_j)$ with $q_j<S_j$, $U_i^+=(C_i,S_i)$, and $S_i(0^+,t)\geq S_i$. Otherwise, $U_i(0^+,t)=U_i^+=(q_i,C_i)$ with $q_i<S_i$, and
\bqs
q_i&=&\frac{D_0}{C_i+S_j(0^+,t)} C_i \geq \frac{C_i}{C_1+C_2}D_0 \geq S_i,
\eqs
which is impossible. Since at least one of the downstream links has SUC stationary state, the stationary states on links $i$ and $j$ are OC and SUC respectively.
 From \refe{evacuationdiverge}, we have a unique interior state on link $i$, $U_i(0^+,t)=(C_i,\frac{S_i}{D_0-S_i}C_j)$, and $q_j=D_0-S_i$.
\ei
\ei
For the four cases, it is straightforward to show that \refe{evacuationdiverge2} always holds.
\eop

\pdfbookmark[1]{Appendix D}{appendixd}
\section*{Appendix D: Proof of Theorem \ref{thm:priorityevacuation}}
{\em Proof}. 

First \refe{priorityevacuation} implies that
\bqs
q_0&=&\min\{S_1(0^+,t)+S_2(0^+,t), D_0(0^-,t) \},
\eqs
which can be shown for three cases: (i) $S_1(0^+,t)+S_2(0^+,t)<D_0(0^-,t)$, (ii) $S_i(0^+,t)\geq \a_iS_0(0^-,t)$, and (iii) $S_1(0^+,t)+S_2(0^+,t)\geq D_0(0^-,t)$ and $S_i(0^+,t)\leq \a_iS_0(0^-,t)$.
\bi
\item [(1)] When $S_1+S_2<D_0$, $q_0=q_1+q_2\leq S_1+S_2<D_0\leq C_0$. Thus the downstream stationary state is SOC with $U_0^-=U_0(0^-,t)=(C_0,q_0)$. In the following, we prove that $q_i=S_i$, which is consistent with \refe{priorityevacuation2}.
\bi
\item[(i)] Assuming that $q_i<S_i\leq C_i$, then the stationary state on link $i$ is SUC with $U_i^+=U_i(0^+,t)=(q_i,C_i)$. From \refe{priorityevacuation}, we have
\bqs
q_i&=&\min\{C_i, \max\{C_0-S_j(0^+,t), \a_i C_0\}\}=\max\{C_0-S_j(0^+,t), \a_i C_0\}<C_i,\\
q_j&=&\min\{S_j(0^+,t), \max\{C_0-C_i, \alpha_j C_0\}\}.
\eqs
We show that the two equations have no solutions for either $\a_j C_0 \leq S_j(0^+,t)$ or $\a_j C_0 > S_j(0^+,t)$. Thus $q_i=S_i$.

\bi
\item [(a)] When $\a_j C_0 \leq S_j(0^+,t)$, we have $\a_i C_0 \geq C_0-S_j(0^+,t)$. From the first equation we have $q_i=\a_i C_0$. From the second equation we have $q_j= S_j(0^+,t)\geq \a_j C_0$ or $q_j=\max\{C_0-C_i, \alpha_j C_0\}\geq \a_j C_0$. Thus $q_i+q_j\geq C_0\geq D_0$, which contradicts $q_0<D_0$.
\item [(b)] When $\a_j C_0 > S_j(0^+,t)$, we have $\a_i C_0 < C_0-S_j(0^+,t)$. From the first equation we have $q_i= C_0-S_j(0^+,t)$. From the second equation we have $q_j=S_j(0^+,t)$. Thus $q_i+q_j=C_0$, which contradicts $q_0<D_0$.
\ei
\ei

\item [(2)] When $S_i \geq \a_i D_0$, $D_0-S_j \leq \a_i D_0$. In the following we show that $q_0=D_0$ and $q_i=\a_i D_0$, which is consistent with \refe{priorityevacuation2}.
\bi
\item [(i)] If $q_0<D_0$, then the stationary state on link 0 is SOC with $U_0^-=U_0(0^-,t)=(C_0,q_0)$. Also at least one of the downstream stationary states is SUC, since, otherwise, $q_1+q_2=S_1+S_2\geq D_0$. Here we assume $U_i^+=U_i(0^+,t)=(q_i,C_i)$.
From  \refe{priorityevacuation} we have
\bqs
q_i&=&\min\{C_i, \max\{C_0-S_j(0^+,t), \alpha_i C_0\}\}=\max\{C_0-S_j(0^+,t), \alpha_i C_0\},\\
q_j&=&\min\{S_j(0^+,t), \max\{C_0-C_i, \alpha_j C_0\}\}.
\eqs 
We show that the two equations have no solutions for either $C_0-S_j(0^+,t) \geq \a_i C_0$ or $C_0-S_j(0^+,t) < \a_i C_0$. Thus $q_0=D_0$.
\bi
\item [(a)] If $C_0-S_j(0^+,t) \geq \a_i C_0$, $S_j(0^+,t) \leq \a_j C_0$. From the first equation we have $q_i= C_0-S_j(0^+,t)$. From the second equation we have $q_j=S_j(0^+,t)$. Thus $q_i+q_j=C_0$, which contradicts $q_0<D_0\leq C_0$.
\item [(b)] If $C_0-S_j(0^+,t) < \a_i C_0$, $S_j(0^+,t) > \a_j C_0$. From the first equation we have $q_i= \a_i C_0$. From the second equation we have $q_j=S_j(0^+,t) > \a_j C_0$ or $q_j=\max\{C_0-C_i, \alpha_j C_0\}\geq \a_j C_0$. Thus $q_i+q_j\geq C_0$, which contradicts $q_0<D_0\leq C_0$.
\ei
\item [(ii)] If $q_i<\a_i D_0\leq S_i\leq C_i$ for any $i=1,2$, then $U_i^+=U_i(0^+,t)=(q_i,C_i)$. From  \refe{priorityevacuation} we have
\bqs
q_i&=&\max\{D_0(0^-,t)-S_j(0^+,t), \alpha_i D_0(0^-,t)\}<C_i,\\
q_j&=&\min\{S_j(0^+,t), \max\{D_0(0^-,t)-C_i, \alpha_j D_0(0^-,t)\}\}.
\eqs
The first equation implies that $\alpha_i D_0(0^-,t)< \a_i D_0$; i.e., $D_0(0^-,t)<D_0$. In addition, $D_0(0^-,t)-S_j(0^+,t)<\a_i D_0$. Thus, $D_0(0^-,t)-C_i<D_0-C_i<D_0-\a_i D_0=\a_j D_0$, and $\max\{D_0(0^-,t)-C_i, \alpha_j D_0(0^-,t)\}<\a_j D_0$. From the second equation we have $q_j<\a_j D_0$. Thus $q_i+q_j<D_0$, which contradicts $q_i+q_j=D_0$. Thus $q_i\geq \a_i D_0$ for $i=1,2$. Since $q_i+q_j=D_0$, $q_i=\a_i D_0$.
\ei

\item[(3)] When $S_i+S_j \geq D_0$ and $S_i\leq \a_i D_0$ for $i,j=1$ or 2 and $i\neq j$. In the following we show that $q_0=D_0$ and $q_i=S_i$, which is consistent with \refe{priorityevacuation2}.
\bi
\item [(i)] If $q_0<D_0$, then the stationary state on link 0 is SOC with $U_0^-=U_0(0^-,t)=(C_0,q_0)$. We first prove that at least one downstream stationary state is SUC and then that none of the downstream stationary states can be SUC. Therefore, $q_0=D_0$.

\bi
\item [(a)] If none of the downstream stationary states are SUC, then $q_1+q_2=S_1+S_2\geq D_0$, which contradicts $q_0<D_0$. Thus, at least one of the downstream stationary states is SUC.
\item [(b)] Assuming that $q_i<S_i$, then  $U_i^+=U_i(0^+,t)=(q_i,C_i)$.
From  \refe{priorityevacuation} we have 
\bqs
q_i&=&\min\{C_i, \max\{C_0-S_j(0^+,t), \alpha_i C_0\}\}=\max\{C_0-S_j(0^+,t), \alpha_i C_0\}<S_i,
\eqs 
which is not possible, since $S_i\leq \a_i D_0$. Thus $q_i=S_i$.
\item [(c)] Assuming that $q_j<S_j$, then  $U_j^+=U_j(0^+,t)=(q_j,C_j)$. Since $q_0=\min\{S_i(0^+,t)+C_j, C_0 \}<D_0$, we have $q_0=S_i(0^+,t)+C_j<D_0$. 
From  \refe{priorityevacuation} we have $S_i=q_i\leq S_i(0^+,t)$. Thus $S_i+S_j\leq S_i(0^+,t)+C_j <D_0$, which contradicts $S_i+S_j\geq D_0$.
\ei

\item [(ii)] If $q_i<S_i$, then $U_i^+=U_i(0^+,t)=(q_i,C_i)$. From \refe{priorityevacuation}, we have 
\bqs
q_i&=&\max\{D_0(0^-,t)-S_j(0^+,t), \alpha_i D_0(0^-,t)\}<S_i\leq \a_i D_0,\\
q_j&=&\min\{S_j(0^+,t), \max\{ D_0(0^-,t)-C_i, \a_j D_0(0^-,t)\} \}.
\eqs
From the first equation we have that $D_0(0^-,t)<D_0$. We show that the two equations have no solutions for either $D_0(0^-,t)-S_j(0^+,t) \geq \alpha_i D_0(0^-,t)$ or $D_0(0^-,t)-S_j(0^+,t) < \alpha_i D_0(0^-,t)$. Therefore $q_i=S_i$.
\bi
\item [(a)] When $D_0(0^-,t)-S_j(0^+,t) \geq \alpha_i D_0(0^-,t)$, we have $S_j(0^+,t) \leq \a_j D_0(0^-,t)$. Thus $q_i=D_0(0^-,t)-S_j(0^+,t)$ and  $q_j=S_j(0^+,t)$. Then $q_i+q_j=D_0(0^-,t)<D_0$, which contradicts $q_i+q_j=D_0$.
\item [(b)] When $D_0(0^-,t)-S_j(0^+,t) < \alpha_i D_0(0^-,t)$, we have $q_i=\a_i D_0(0^-,t)$ and $D_0(0^-,t)-C_i < D_0(0^-,t)-q_i=\a_j D_0(0^-,t)$. Thus $q_j\leq \a_j D_0(0^-,t)$, and $q_i+q_j\leq D_0(0^-,t)<D_0$, which contradicts $q_0=D_0$. 
\ei

\ei

\ei

\eop

\begin{thebibliography}{}

\bibitem[Bultelle et~al., 1998]{bultelle1998shock}
Bultelle, M., Grassin, M., and Serre, D. (1998).
\newblock {Unstable Godunov discrete profiles for steady shock waves}.
\newblock {\em SIAM Journal on Numerical Analysis}, 35(6):2272--2297.

\bibitem[Cassidy, 2003]{cassidy2003freeway}
Cassidy, M. (2003).
\newblock {Freeway On-Ramp Metering, Delay Savings, and Diverge Bottleneck}.
\newblock {\em Transportation Research Record: Journal of the Transportation
  Research Board}, 1856:1--5.

\bibitem[Coclite et~al., 2005]{coclite2005network}
Coclite, G., Garavello, M., and Piccoli, B. (2005).
\newblock {Traffic flow on a road network}.
\newblock {\em SIAM Journal on Mathematical Analysis}, 36:1862.

\bibitem[Colella and Puckett, 2004]{collela2004_pde}
Colella, P. and Puckett, E.~G. (2004).
\newblock {\em Modern Numerical Methods for Fluid Flow}.
\newblock In draft.

\bibitem[Courant et~al., 1928]{courant1928CFL}
Courant, R., Friedrichs, K., and Lewy, H. (1928).
\newblock {{\\"U}ber die partiellen Differenzengleichungen der mathematischen
  Physik}.
\newblock {\em Mathematische Annalen}, 100:32--74.

\bibitem[Daganzo, 1999]{daganzo1999remarks}
Daganzo, C. (1999).
\newblock {Remarks on Traffic Flow Modeling and Its Appi}.
\newblock {\em Traffic and mobility: simulation, economics, environment}, page
  105.

\bibitem[Daganzo, 1994]{daganzo1994ctm}
Daganzo, C.~F. (1994).
\newblock The cell transmission model: a dynamic representation of highway
  traffic consistent with hydrodynamic theory.
\newblock {\em Transportation Research Part B}, 28(4):269--287.

\bibitem[Daganzo, 1995]{daganzo1995ctm}
Daganzo, C.~F. (1995).
\newblock The cell transmission model \m{II}: Network traffic.
\newblock {\em Transportation Research Part B}, 29(2):79--93.

\bibitem[Daganzo, 1997]{daganzo1997special}
Daganzo, C.~F. (1997).
\newblock A continuum theory of traffic dynamics for freeways with special
  lanes.
\newblock {\em Transportation Research Part B}, 31(2):83--102.

\bibitem[Daganzo et~al., 1999]{daganzo1999phase}
Daganzo, C.~F., Cassidy, M.~J., and Bertini, R.~L. (1999).
\newblock Possible explanations of phase transitions in highway traffic.
\newblock {\em Transportation Research A}, 33:365--379.

\bibitem[Daganzo et~al., 1997]{daganzo1997it}
Daganzo, C.~F., Lin, W.-H., and {Del Castillo}, J.~M. (1997).
\newblock A simple physical principle for the simulation of freeways with
  special lanes and priority vehicles.
\newblock {\em Transportation Research Part B}, 31(2):103--125.

\bibitem[{Del Castillo} and Benitez, 1995]{delcastillo1995fd_empirical}
{Del Castillo}, J.~M. and Benitez, F.~G. (1995).
\newblock On the functional form of the speed-density relationship - \m{II}:
  Empirical investigation.
\newblock {\em Transportation Research Part B}, 29(5):391--406.

\bibitem[Engquist and Osher, 1980]{engquist1980calculation}
Engquist, B. and Osher, S. (1980).
\newblock {Stable and entropy satisfying approximations for transonic flow
  calculations}.
\newblock {\em Mathematics of Computation}, 34(149):45--75.

\bibitem[Fazio et~al., 1990]{fazio1990behavioral}
Fazio, J., Michaels, R., Reilly, W., Schoen, J., and Poulis, A. (1990).
\newblock {Behavioral model of freeway exiting}.
\newblock {\em Transportation Research Record}, 1281:16--27.

\bibitem[{Federal Highway Administration}, 2004]{fhwa2004tft}
{Federal Highway Administration} (2004).
\newblock {\em {Traffic Flow Theory: A State of the Art Report}}.
\newblock Transportation Research Board.

\bibitem[Greenshields, 1935]{greenshields1935capacity}
Greenshields, B.~D. (1935).
\newblock A study in highway capacity.
\newblock {\em Highway Research Board Proceedings}, 14:448--477.

\bibitem[Haberman, 1977]{haberman1977model}
Haberman, R. (1977).
\newblock {\em Mathematical models}.
\newblock Prentice Hall, Englewood Cliffs, NJ.

\bibitem[Holden and Risebro, 1995]{holden1995unidirection}
Holden, H. and Risebro, N.~H. (1995).
\newblock A mathematical model of traffic flow on a network of unidirectional
  roads.
\newblock {\em SIAM Journal on Mathematical Analysis}, 26(4):999--1017.

\bibitem[Jin, 2010]{jin2010_merge}
Jin, W.-L. (2010).
\newblock {Continuous kinematic wave models of merging traffic flow}.
\newblock {\em Transportation Research Part B}.
\newblock In Press.

\bibitem[Jin et~al., 2009]{jin2009sd}
Jin, W.-L., Chen, L., and Puckett, E.~G. (2009).
\newblock {Supply-demand diagrams and a new framework for analyzing the
  inhomogeneous Lighthill-Whitham-Richards model}.
\newblock {\em {Proceedings of the 18th International Symposium on
  Transportation and Traffic Theory}}, pages 603--635.

\bibitem[Jin and Zhang, 2003a]{jin2003inhlwr}
Jin, W.-L. and Zhang, H.~M. (2003a).
\newblock The inhomogeneous kinematic wave traffic flow model as a resonant
  nonlinear system.
\newblock {\em Transportation Science}, 37(3):294--311.

\bibitem[Jin and Zhang, 2003b]{jin2003merge}
Jin, W.-L. and Zhang, H.~M. (2003b).
\newblock On the distribution schemes for determining flows through a merge.
\newblock {\em Transportation Research Part B}, 37(6):521--540.

\bibitem[Jin and Zhang, 2004]{jin2004network}
Jin, W.-L. and Zhang, H.~M. (2004).
\newblock A multicommodity kinematic wave simulation model of network traffic
  flow.
\newblock {\em Transportation Research Record: Journal of the Transportation
  Research Board}, 1883:59--67.

\bibitem[Lebacque and Khoshyaran, 2005]{lebacque2005network}
Lebacque, J. and Khoshyaran, M. (2005).
\newblock {First order macroscopic traffic flow models: Intersection modeling,
  Network modeling}.
\newblock {\em Proceedings of the 16th International Symposium on
  Transportation and Traffic Theory}, pages 365--386.

\bibitem[Lebacque, 1996]{lebacque1996godunov}
Lebacque, J.~P. (1996).
\newblock {The Godunov scheme and what it means for first order traffic flow
  models}.
\newblock {\em Proceedings of the 13th International Symposium on
  Transportation and Traffic Theory}, pages 647--678.

\bibitem[Lighthill and Whitham, 1955]{lighthill1955lwr}
Lighthill, M.~J. and Whitham, G.~B. (1955).
\newblock On kinematic waves: \m{II. A} theory of traffic flow on long crowded
  roads.
\newblock {\em Proceedings of the Royal Society of London A},
  229(1178):317--345.

\bibitem[Liu et~al., 1996]{liu1996junction}
Liu, G., Lyrintzis, A.~S., and Michalopoulos, P.~G. (1996).
\newblock Modelling of freeway merging and diverging flow dynamics.
\newblock {\em Applied Mathematical Modelling}, 20.

\bibitem[Mu{\~ n}oz and Daganzo, 2002]{munoz2002diverge}
Mu{\~ n}oz, J.~C. and Daganzo, C.~F. (2002).
\newblock The bottleneck mechanism of a freeway diverge.
\newblock {\em Transportation Research Part A}, 36(6):483--505.

\bibitem[Newell, 1993]{newell1993sim}
Newell, G.~F. (1993).
\newblock A simplified theory of kinematic waves in highway traffic \m{I}:
  General theory. \m{II}: Queuing at freeway bottlenecks. \m{III}:
  Multi-destination flows.
\newblock {\em Transportation Research Part B}, 27(4):281--313.

\bibitem[Newell, 1999]{newell1999exit}
Newell, G.~F. (1999).
\newblock Delays caused by a queue at a freeway exit ramp.
\newblock {\em Transportation Research Part B}, 33:337--350.

\bibitem[Ngoduy et~al., 2006]{ngoduy2006continuum}
Ngoduy, D., Hoogendoorn, S., and Van~Zuylen, H. (2006).
\newblock {Continuum traffic model for freeway with on-and off-ramp to explain
  different traffic-congested states}.
\newblock {\em Transportation research record}, (1965):91--102.

\bibitem[Papageorgiou, 1990]{papageorgiou1990assignment}
Papageorgiou, M. (1990).
\newblock Dynamic modelling, assignment and route guidance in traffic networks.
\newblock {\em Transportation Research Part B}, 24(6):471--495.

\bibitem[Qiu and Jin, 2008]{qiu2008see}
Qiu, K.-F. and Jin, W.-L. (2008).
\newblock {Studies of Emergency Evacuation Strategies based on Kinematic Wave
  Models of Network Vehicular Traffic}.
\newblock In {\em Intelligent Transportation Systems, 2008. ITSC 2008. 11th
  International IEEE Conference on}, pages 222--227.

\bibitem[Richards, 1956]{richards1956lwr}
Richards, P.~I. (1956).
\newblock Shock waves on the highway.
\newblock {\em Operations Research}, 4(1):42--51.

\bibitem[Sheffi et~al., 1982]{Sheffi1982evacuation}
Sheffi, Y., Mahmassani, H., and Powell, W. (1982).
\newblock {A transportation network evacuation model}.
\newblock {\em Transportation Research Part A}, 16(3):209--218.

\bibitem[van Leer, 1984]{vanleer1984upwind}
van Leer, B. (1984).
\newblock {On the relation between the upwind-differencing schemes of Godunov,
  Engquist-Osher and Roe}.
\newblock {\em SIAM Journal on Scientific and Statistical Computing},
  5(1):1--20.

\end{thebibliography}
\end {document}